\documentclass[reqno,11pt]{amsart}
\usepackage{hyperref}
\usepackage{geometry}
\usepackage{enumitem}
\usepackage{amsmath}
\usepackage{amssymb}
\usepackage{amsthm}
\usepackage{amsrefs}
\usepackage{amsfonts}
\usepackage{xcolor}

\usepackage{mathtools}

\usepackage{tikz}
\usepackage{stackengine}

\def\XXint#1#2#3{{\setbox0=\hbox{$#1{#2#3}{\int}$ }
\vcenter{\hbox{$#2#3$ }}\kern-.6\wd0}}

\makeatletter
\newcommand*{\rom}[1]{\expandafter\@slowromancap\romannumeral #1@}

\newcommand{\di}{\text{d}}
\newcommand{\diam}{\text{diam}}
\newcommand{\dist}{\text{dist}}
\newcommand{\Int}{\text{Int}}

\newcommand{\PU}{\text{PU}}
\newcommand{\R}{\mathbb{R}}

\newcommand{\Z}{\mathbb{Z}}
\newcommand{\Sphere}{\mathbb{S}}
\newcommand{\h}{\mathbb{H}}

\newcommand{\N}{\mathbb{N}}

\newcommand{\sm}{\setminus}

\newcommand{\bthm}{\begin{thm}}
\newcommand{\ethm}{\end{thm}}
\newcommand{\bproof}{\begin{proof}}
\newcommand{\eproof}{\end{proof}}
\newcommand{\blem}{\begin{lem}}
\newcommand{\elem}{\end{lem}}
\newcommand{\brem}{\begin{rem}}
\newcommand{\erem}{\end{rem}}
\newcommand{\eeqn}{\end{equation}}
\newcommand{\eeqnn}{\end{equation*}}
\newcommand{\beqn}{\begin{equation}}
\newcommand{\beqnn}{\begin{equation*}}
\newcommand{\eprop}{\end{prop}}
\newcommand{\eexm}{\end{exm}}
\newcommand{\enexm}{\end{nexm}}
\newcommand{\ecor}{\end{cor}}
\newcommand{\bcor}{\begin{cor}}
\newcommand{\bexm}{\begin{exm}}
\newcommand{\bnexm}{\begin{nexm}}
\newcommand{\bprop}{\begin{prop}}
\newcommand{\bdefn}{\begin{defn}}
\newcommand{\edefn}{\end{defn}}
\newcommand{\benum}{\begin{enumerate}}
\newcommand{\eenum}{\end{enumerate}}

\title[Diophantine Approximation and Large Intersections]{Diophantine Approximation, large intersections and geodesics in negative curvature}

\begin{document}
\theoremstyle{plain}
\newtheorem{thm}{Theorem}[section]
\newtheorem{lem}[thm]{Lemma}
\newtheorem{prop}[thm]{Proposition}
\newtheorem{cor}[thm]{Corollary}

\theoremstyle{definition}
\newtheorem{defn}[thm]{Definition}
\newtheorem{exm}[thm]{Example}
\newtheorem{nexm}[thm]{Non Example}
\newtheorem{prob}[thm]{Problem}

\theoremstyle{remark}
\newtheorem{rem}[thm]{Remark}

\author{Anish Ghosh}
\address{\textbf{Anish Ghosh} \\
School of Mathematics,
Tata Institute of Fundamental Research, Mumbai, India 400005}
\email{ghosh@math.tifr.res.in}

\author{Debanjan Nandi}
\address{\textbf{Debanjan Nandi} \\
The Weizmann Institute of Science, 234 Herzl Street, Rehovot 76100, Israel}
\email{debanjan.nandi@weizmann.ac.il}
\date{}

\thanks{A.\ G.\ gratefully acknowledges support from a grant from the Indo-French Centre for the Promotion of Advanced Research; and a MATRICS grant from the Science and Engineering Research Board and a grant from the Infosys foundation. Both authors gratefully acknowledge support from a Department of Science and Technology, Government of India, Swarnajayanti fellowship}

\date{}

\begin{abstract}
In this paper we prove quantitative results about geodesic approximations to submanifolds in negatively curved spaces.  
Among the main tools is a new and general Jarn\'{i}k-Besicovitch type theorem in Diophantine approximation. The framework we develop is flexible enough to treat manifolds of variable negative curvature, a variety of geometric targets, and logarithm laws as well as spiraling phenomena in both measure and dimension aspect. Several of the results are new also for manifolds of constant negative sectional curvature. We further establish a large intersection property of Falconer in this context.

\end{abstract}

\maketitle
\tableofcontents

\section{Introduction}

This paper has its origins in fundamental work of Patterson \cite{Pat2} on Diophantine approximation and Fuchsian groups, and of Sullivan \cite{Sul} on logarithm laws for cuspidal excursions of geodesics in hyperbolic manifolds. Our subject here is a quantitative study of certain asymptotic properties of geodesics in negatively curved spaces. 
Consider a negatively curved manifold $M$ and a `target' $N$ which is a subset of the manifold. Given a point $p$ in $M$, what can be said about the size (in terms of measure and dimension) of the set of geodesics starting at $p$ which
\benum
\item infinitely often visit neighbourhoods of $N$ which are shrinking in volume? We refer to this problem as the \emph{shrinking target problem}, following Hill and Velani \cite{HV} as well as Kleinbock and Margulis \cite{KMinv}.
\item take longer and longer sojourns into a fixed neighbourhood of $N$? We refer to this as the \emph{spiral trap problem}.
\eenum
When $N$ is an `isolated point at infinity', (1) is the setting of Sullivan's celebrated work \cite{Sul}. A Hausdorff dimension estimate in this case is given in the work of Melian-Pestana \cite{MP}. When $N$ is a point in $M$, (1) has been considered in the work of Maucourant \cite{Mau}, who provides a zero-one law (i.e. a measure theoretic Borel-Cantelli statement) in the case $M$ is a finite volume manifold of constant negative sectional curvature. When $N$ is a closed geodesic bounding a funnel in a surface of constant negative sectional curvature without cusps, a Hausdorff dimension estimate for (1) is provided by Dodson, Meli\'{a}n, Pestana and Velani \cite{DMPV}. When $M$ is a closed manifold of variable, strictly negative sectional curvature, and $N$ is again a point in $M$, a zero-one law for (1) is known from the work of Hersonsky, Paulin and Aravinda in \cite{HP2} and a Hausdorff dimension estimate is known from the work of Hersonsky-Paulin \cite{HP1}. 
Question (2) has been considered in the work of Hersonsky-Paulin \cite{HP2} who give a zero-one law for convex subsets in the more general setting of CAT($-1$) spaces. A necessarily incomplete list of work related to the stated problems includes the papers \cite{Strat, Velani, DMPV, BDV, AGP, AP, PP}.\\

In this paper, we develop a framework for obtaining comprehensive measure and dimension results for the shrinking target problem and the spiral trap problem. 
This allow us to treat general and natural geometric targets and our results are valid in a very general setting, allowing for instance, manifolds of variable negative curvature and more generally, CAT($-1$) metric spaces. Nevertheless, several of our results are new even for closed manifolds of constant negative curvature. Namely, we give sharp Hausdorff dimension estimates for (2). For (1), we relax the requirement that the manifold be closed, and generalize from targets being points to more general convex subsets, and obtain zero-one laws as well as sharp dimension estimates. Let us detail and further illustrate the type of problems studied here with the following very special cases of results proved in \S \ref{sec:spiraling}. For a smooth manifold $M$ and $x\in M$, let $SM_x$ denote the unit tangent sphere at $x$.

Theorem 1.1 below is a special case of Theorem \ref{loglaw1}.
\bthm \label{loglaw2} Let $M$ be a closed manifold of dimension $n$ with constant negative sectional curvature, $k=-1$. Let $N$ be a compact, totally geodesic submanifold of $M$ of dimension $0\leq s\leq n-1$. Let $\tau\geq 0$ be fixed. Then given $x_0\in M$, we have that the set 
\[
  E^{\tau}_{N} =\left\lbrace v\in SM_{x_0} \;\middle|\;
  \begin{tabular}{@{}l@{}}
    $\exists\; \text{positive times}\; t_n\to\infty\;\text{such that}$\\ $\gamma_v(t_n)\subset B(N,Ce^{-\tau t_n})\;\text{for some} \; C>0$
   \end{tabular}
  \right\rbrace
\] has Hausdorff dimension $$ \dim_{\mathcal{H}}(E_{N}^{\tau})=(n-1)\cdot\frac{1+\frac{\tau s}{n-1}}{1+\tau},$$ where $\gamma_v$ is the geodesic at $x_0$ at time zero with direction $v$.
\ethm
Theorem 1.2 is a special case of Theorem \ref{loglaw2gen}.
\bthm Let $\epsilon>0$ and $\tau\geq 0$ be fixed.
Let $M$ be a closed manifold of dimension $n$ with constant negative sectional curvature, $k= -1$. Let $N$ be a compact, totally geodesic submanifold of dimension $1\leq s\leq n-1$. Let $x_0\in M$. Then the set 
 \[
  E^{\tau}_{N} =\left\lbrace v\in SM_{x_0} \;\middle|\;
  \begin{tabular}{@{}l@{}}
    $\exists\; \text{positive times}\; t_n \to\infty\;\text{such that}$\\ $\gamma_v(t_n,t_n+\tau t_n)\subset B(N,\epsilon)$
   \end{tabular}
  \right\rbrace
\] has $$ \text{dim}_{\mathcal{H}} (E^{\tau}_{N})=(n-1)\cdot \frac{1+\frac{\tau(s-1)}{n-1}}{1+\tau},$$ where $\gamma_v$ is the geodesic at $x_0$ at time zero with direction $v$.
 \ethm

Theorem 1.3 is a special case of Corollary \ref{full_dim}.

\bthm
Let $M$ be a closed manifold of dimension $n$ with constant negative curvature $k=-1$. Let $N$ be a compact, totally geodesic submanifold in $M$ of dimension $0\leq s \leq n-1$. Let $x_0\in M$. Fix $\epsilon>0$.
\benum 
\item  If $1\leq s\leq n-1$, then
 \[
  E_{N} =\left\lbrace v\in SM_{x_0} \;\middle|\;
  \begin{tabular}{@{}l@{}}
    $\exists\; \text{positive times}\; t_n \to\infty\;\text{such that}$\\ $\gamma_v(t_n,t_n+\tau t_n)\subset B(N,\epsilon)\;\text{for some}\;\tau>0$
   \end{tabular}
  \right\rbrace
\] has $$ \dim_{\mathcal{H}} (E_{N})= n-1.$$ 
\item Also, for $0\leq s\leq n-1$
\[
  E'_{N} =\left\lbrace v\in SM_{x_0} \;\middle|\;
  \begin{tabular}{@{}l@{}}
    $\exists\; \text{positive times}\; t_n\to\infty\;\text{such that}$\\ $\gamma_v(t_n)\subset B(N,e^{-\tau t_n})\;\text{for some} \; \tau>0$
   \end{tabular}
  \right\rbrace
\] has  $$ \dim_{\mathcal{H}}(E'_{N})= n-1,$$ 
\eenum
where $\gamma_v$ is the geodesic at $x_0$ at time zero with direction $v$. The set $E_N$ is a null set for the Lebesgue measure in $SM_{x_0}$ for $1\leq s \leq n-1$. The set $E_N'$ is a null set for the Lebesgue measure if and only if $0\leq s<n-1$. 
\ethm

\subsection{A correspondence between dynamics of the geodesic flow and Diophantine approximation}
Asymptotic properties of geodesics in hyperbolic manifolds, namely spiral traps and shrinking target properties, are closely related to metric Diophantine approximation. 

Namely, these phenomena can be recast as Diophantine approximation questions about (subsets of) orbits of the fundamental group in the visual boundary of the universal cover of the manifold. For instance, Khintchine type results in Diophantine approximation follow as a result of mixing of the geodesic flow. On the other hand, Jarnik-Besicovitch type results utilizing the metric measure structure of the limit set can be used to obtain quantitative estimates on fine asymptotic behaviour of geodesics. There is thus a bilateral correspondence between these two rich areas.
In this paper we prove some of the most general results exploring this correspondence and obtain sharp estimates for the problems stated above. 

We now describe this connection in greater detail. Recall that the Hausdorff measure and dimension refinements of Khintchine's theorem in Diophantine approximation were obtained by Jarn\'{i}k \cite{Jar}, and independently by Besicovitch \cite{Bes}. In more recent fundamental work, Beresnevich and Velani \cite{BV} proved a \emph{mass transference principle} providing a surprising connection between Lebesgue measure theoretic statements (like Khintchine's theorem) about certain `limsup' sets and Hausdorff measure theoretic statements (like the Jarn\'{i}k-Besicovitch theorem) for some null subsets. Loosely speaking, the mass transference principle is the assertion that the former implies the latter. 
In section 3, we prove a general mass transference principle, Theorem \ref{density} which is well suited  to geometric considerations. In particular, it can be used to  study the action of a discrete group $\Gamma$ acting properly on a CAT(-1) space $X$ and ergodically on its limit set in the visual boundary. We 
then connect metric aspects of the $\Gamma$ action on the visual boundary of $X$ and the asymptotic behaviour of geodesics on $X/\Gamma$. Our result may also be viewed as a very general version of the Jarn\'{i}k-Besicovitch theorem. 
In contrast to the aforementioned works, we do not fix the shape of the sets determining the exceptional sets under consideration; our mass transference principle applies to the case of open sets with \emph{arbitrary} shapes, as necessitated by our setting.

\subsection{Large Intersections}
The sets $E^{\tau}_{N}$ of geodesic directions hitting exponentially shrinking targets in Theorem \ref{loglaw2} above, are zero measure sets and so there is no control \emph{a priori} on the size of their intersections which could be trivial. Nevertheless, they have a `large intersection' property which leads to the surprising fact that countable intersections have Hausdorff dimensions bounded below by the infimum of the respective Hausdorff dimensions of the individual sets (see Theorem \ref{LIP_intro} below).

In \cite{Fal}, Falconer defined a class of subsets, denoted $\mathcal{G}^s,$ of $\mathbb{R}^n,$  which form a  maximal class of $G_{\delta}$-sets of dimension at least $s$ that is closed under countable intersections and under similarity transformations (see also \cite{Fal85}). He named this property the large intersection property. Falconer's definition unifies several earlier categories of sets with similar properties including the `regular systems' of Baker and Schmidt \cite{BakSch} and the `ubiquitous systems' of Dodson, Rynne and Vickers \cite{DRV} and consequently these classes play an important role in Diophantine approximation. Theorem \ref{LIP} provides in particular a large intersection property for actions of hyperbolic groups on the visual boundaries of hyperbolic metric spaces. We avoid the use of net measures as we use the `coarser' dyadic decomposition coming from shadows of balls centered at orbits of $\Gamma$; instead we use a variant of the Hausdorff measure. Here, by coarse we mean that the decomposition does not need to have the property that when sets of different `generations' overlap, then the interior of one has to be completely contained in the other. We note however that more abstract versions (not determined by concrete geometric properties) of dyadic decompositions also exist for Ahlfors-regular metric spaces, see \cite{HytKar}, \cite{Chr}. The abstract dyadic decomposition given in \cite{HytKar}, also has the properties we require. The following result is an application of Theorem \ref{LIP} to the problem of spiraling of geodesics. It is a special case of Theorem \ref{spiral_lip}.  

\bthm\label{LIP_intro}
Let $M$ be a closed manifold of dimension $n$, with pinched sectional curvature $-a^2\leq k\leq -1$. Let $N_i$ be a countable collection of closed totally geodesic submanifolds or points of $M$. Let $x_0\in M$. Let $\tau_i$ be a sequence of positive numbers. Then the set of directions \[
  E =\left\lbrace v\in SM_{x_0} \;\middle|\;
  \begin{tabular}{@{}l@{}}
    $\exists\; \text{positive times}\; t_n^{(i)} \to\infty\;\text{for each $i$ such that}$\\ $\gamma_v(t_n^{(i)})\in B(N_i,e^{-\tau_i t_n^{(i)}})\;\text{for each i}$
   \end{tabular}
  \right\rbrace
\] has Hausdorff dimension $$ \inf_i\;(n-1)\frac{1+\frac{\tau_i\dim(N_i)}{n-1}}{1+\tau_i}\leq \dim_{\mathcal{H}}(E)\leq \inf_i\;(n-1)\frac{1+\frac{\tau_i\dim(N_i)}{n-1}}{1+\tau_i/a}.$$
\ethm
Note that the dimension lower bound above is always positive as long as either the sequence $\tau_i$ is bounded above or none of the submanifolds $N_i$ are singletons.\\

Finally, we note that a detailed and systematic study of Diophantine approximation in the context of hyperbolic groups has been carried out in the monograph \cite{FSU} of Fishman, Simmons and Urbanski. They construct `partition structures' (similar to our dyadic decomposition in the case of geometric actions) and prove results about the limit set of discrete group actions in very general settings. The dyadic decomposition we use enjoys stronger properties that we require for establishing the large intersection property. The aforementioned paper has no intersection with our results.

\subsection*{Structure of the paper} In the next section, we set some notation, gather metric and measure preliminaries and introduce the Whitney decomposition for metric spaces with a dyadic decomposition. In section $3$, we  introduce an abstract framework for Diophantine approximation and study it using adjusted well distributed systems, which we also introduce. Section $4$ is devoted to results about Falconer's large intersection property and section $5$ is devoted to our results on spiraling of geodesics and associated $0-1$ laws. 
Finally, we will use the notation $A \gtrsim B$ to mean that there is a constant $c > 0$ such that $A \geq cB$. Dependence of the constant on parameters will be specified. The notation $A \approx B $ will stand for $A\gtrsim B \gtrsim A.$

\subsection*{Acknowledgements} We would like to thank Mahan Mj for helpful comments and encouragement. AG would like to thank Yann Bugeaud and Arnaud Durand for helpful conversations about the large intersection property and Fran\c{c}ois Maucourant for answering some questions. We would like to thank Fr\'{e}d\'{e}ric Paulin for his comments on a preliminary version of this paper. Part of this work was completed when both authors were at the International Centre for Theoretical Sciences, Bengaluru as part of the programme \emph{Smooth and Homogeneous Dynamics}. The hospitality of ICTS is gratefully acknowledged.     

\section{Preliminaries}\label{Preliminaries}
Let $(X,\rho)$ be a proper, geodesic and hyperbolic metric space. 

\subsection{The visual metric}
The visual boundary $\partial X$ is a compact metric space with a family of visual metrics $\di_{\beta}$ (mutually \textit{quasisymmetric}), $0<\beta<\beta_X$, for some $\beta_X>0$ which satisfy $$\frac{1}{C_X(\beta)}e^{-\beta(\xi|\eta)}\leq\di_\beta(\xi,\eta)\leq C_X(\beta)e^{-\beta(\xi|\eta)},$$ for $\xi,\eta\in \partial X$, where $$(x|y):=(x|y)_{x_0}=\frac{1}{2}(\rho(x,x_0)+\rho(y,x_0)-\rho(x,y))$$ for $x,y\in X,$ is the Gromov product extended to $\partial X$ by taking limits (to be understood as the distance a pair of geodesic rays joining $x_0$ to the points $x$ and $y$ travel `together'). We fix a $\beta\in (0,\beta_X)$ and let $\di=\di_\beta$ be the corresponding metric in what follows. See Bridson-Haefliger \cite{BriHae}.

\bdefn[Shadow] The set $$S(x,B(z,R)):=\{\xi\in \partial X: \gamma_\xi\cap B(z,R)\neq\emptyset\}$$ is the shadow of the ball $B(z,R)$ with respect to $x$. Here $\gamma_\xi$ represents a geodesic ray from $x$ to $\xi$. 
\edefn

\subsection*{Notation}
Fix $x_0\in X$. We write $$S(g,R) := S(x_0,B(gx_0,R))$$ below.
For $\xi\in\partial X$ we will write $\gamma_\xi$ for a geodesic ray joining $x_0$ to $\xi$. For $x,y$ in $X\cup \partial X$, we write $\gamma_{x,y}$ for a geodesic line (or segment) joining $x$ and $y$. All geodesics (segments, rays and lines) are considered with the unit-speed parametrization.

We now record the following well known fact.
\blem\label{diameter_shadow}
For all $g$ with $|g|_S$ large enough, we have
$$\frac{1}{c_\Gamma}\cdot e^{-\beta\cdot\rho(x_0,gx_0)}\leq \diam(S(g,R))\leq c_\Gamma\cdot e^{-\beta\cdot\rho(x_0,gx_0)}.$$ 
\elem
\bproof
Let $\xi,\xi'\in S(g,R).$ Let $x\in \gamma_\xi\cap B(gx_0,R)$ and $x'\in\gamma_{\xi'}\cap B(gx_0,R)$. Then the claim follows from the inequalities $$(\xi |\xi')\geq \min\{(\xi | x),(x | x'),(x' | \xi')\}-c_\Gamma,$$ and $$(x | x')\geq \min\{(x | \xi),(\xi | \xi'),(\xi' | x')\}-c_\Gamma.$$
\eproof

\subsection{The Patterson-Sullivan (quasiconformal) measure and Ahlfors regularity}
Let $\Gamma$ be a discrete group acting isometrically and properly on $X$. Denote by $\Lambda_\Gamma$ its limit set in $\partial X$. There exist  $\Gamma$-equivariant Patterson-Sullivan (quasiconformal density of)  probability measures $\{\mu_x\}_{x\in X}$, with mutual Radon-Nikodym derivatives given in terms of the Busemann function.
We will only need $\mu$, the Patterson-Sullivan measure in $\Lambda_{\Gamma}$ with respect to our fixed $x_0$. Recall that the critical exponent of $\Gamma$ acting on $X$ is the number (independent of base-point) $$v_{\Gamma} :=\overline{\lim_n}\;\, \frac{1}{n}\cdot \log\,(\#\{g\in\Gamma: gx_0\in B(x_0,n)\}).$$

\blem[Sullivan's Shadow Lemma]\label{shadow_lem}
There exists $R_\Gamma>1$ such that for all $R>R_\Gamma$ and $|g|_S>10\cdot R_\Gamma$, there exists $a_\Gamma\geq 1$ such that 
$$\frac{1}{a_\Gamma}\cdot e^{-v_\Gamma\cdot \rho(x_0,gx_0)}\leq \mu(S(g,R))\leq a_\Gamma\cdot e^{-v_\Gamma\cdot\rho(x_0,gx_0)}.$$
\elem

\blem\label{metric-measure}
There exists $R_\Gamma>1$ such that for all $R>R_\Gamma$ and $|g|_S>R_\Gamma$, there exists $a_\Gamma\geq 1$ such that 
$$\frac{1}{a_\Gamma}\cdot \diam(S(g,R))^{v_\Gamma/\beta}\leq \mu(S(g,R))\leq a_\Gamma\cdot  \diam(S(g,R))^{v_\Gamma/\beta}.$$
\elem
\bproof
This follows from Lemmas \ref{diameter_shadow} and \ref{shadow_lem}.
\eproof
The lemmas above, originally due to Sullivan \cite{Sul1} for $\h^n$, are proved in Coornaert \cite{Coo} in the generality that we consider. The measures were first introduced in the setting of Fuchsian action on the upper-half plane by Patterson \cite{Pat1}.  

\bdefn[Ahlfors regularity]\label{Ahlfors}
We say that a metric measure space $(X,d,\mu)$ is $(c,D)$-Ahlfors regular, if for each $x\in X$ and $0<r<\diam(X)$, $$\frac{1}{c}\cdot r^D\leq \mu(B(x,r))\leq c\cdot r^D.$$
\edefn

\brem\label{regularity_constant}Lemma \ref{metric-measure} says that $(\partial X, \di, \mu)$ is $(c_\Gamma, v_\Gamma/\beta)$-Ahlfors regular.
\erem

\subsection{Dyadic and Whitney decompositions}

The space $(\partial X,\di)$ is locally Ahflors regular. Such spaces admit a decomposition similar to the standard dyadic decomposition of $\R^n$ (see \cite{HytKar} and also \cite{Chr}).

\bdefn[\textbf{Dyadic} decomposition]\label{dyadic}
Let $Y$ be a metric space. A dyadic decomposition $\mathcal{D}$ of $Y$ is a countable collection of subsets $Q_i$, such that there exist constants $A,B,C\geq 1$, $n_0\in\N$, and a decomposition $\mathcal{D}=\bigcup_n\mathcal{W}_n$, where
\benum[label=\text{\ref{dyadic}.\arabic*}]
\item[(2.11.1)] For each $n\in\N$, $\#\mathcal{W}_n<\infty$ and $$\partial X = \underset{\mathcal{W}_n}\bigcup \;Q_i.$$ Moreover, $\lim_n\max\{\diam(Q):Q\in\mathcal{W}_n\}\to 0$.
\item[(2.11.2)] For each $Q_i\in\mathcal{W}$, there exists $x_i\in Q_i$, such that, $$B(x_i,\diam(Q_i)/A)\subset Q_i \subset B(x_i,A\diam(Q_i)).$$
\item[(2.11.3)] For each $n\in\N$, $n\geq n_0$ and $Q_i\in\mathcal{W}_n$, there exists $Q_j\in\mathcal{W}_{n-n_0}$ such that $Q_i\subset Q_j$.
\item[(2.11.4)] Given $n\in\N, l\in\N\cup\{0\}$ and $Q_i\in\mathcal{W}_n$, $$\#\{Q_j\in\mathcal{W}_{n+l}:Q_j\cap Q_i\neq\emptyset\}\leq B^l.$$
\item[(2.11.5)] For $n\in\N$, if $Q\in\mathcal{W}_n$ and $Q'\in\mathcal{W}_{n+1}$, then $$\frac{1}{C}\cdot \diam(Q')\leq\diam(Q)\leq C\cdot \diam(Q').$$
\eenum
\edefn

The following lemma is folklore. The main ingredients are the Milnor-Svar\'c lemma and Sullivan's shadow lemma. The point is that a decomposition in our case is available in terms of concrete objects whose measures and geometry are sufficiently well understood for our purposes.
\blem\label{dyadic_lemma}
Suppose that $X$ is a proper, geodesic, hyperbolic metric space with a geometric action of a group $\Gamma$. Then there exists a dyadic decomposition for the visual boundary $(\partial X,\di)$, equipped with the visual metric given by $$\mathcal{D}=\{S(g,R)\}_{|g|\geq R'},$$ for fixed $x_0\in X$, $R=R(\Gamma, X)>0$ and $R'=R'(\Gamma, X)>0$. The associated constants depend only on $X$ and $\Gamma$.
\elem
A dyadic decomposition has the following property (which follows from definition).
\blem
\label{dyadic_overlap} Given $R\geq 1$ there exists $M=M(R)$, such that for all $Q_i\in \mathcal{D}$, \[\#\left\lbrace Q_j\in\mathcal{D}\middle|\;
\begin{tabular}{@{}l@{}}
$\frac{1}{R}\cdot\diam(Q_i)\leq \diam(Q_j)\leq R\cdot\diam(Q_i)$\\ $\text{and}\;Q_j\cap B(Q_i,R\diam(Q_i))\neq\emptyset$\end{tabular}\right\rbrace \leq M.\]
\elem

Given a dyadic decomposition we also have a \textit{Whitney} decomposition.
\bdefn[\textbf{Whitney} decomposition]\label{Whitney}
Let $Y$ be a bounded metric space and $\mathcal{D}$ be a dyadic decomposition. Let $U\subsetneq Y$ be an open set. A sub-collection $\mathcal{W}$ of $\mathcal{D}$ is a $A$-Whitney decomposition of $U$ if 
\benum[label=\text{\ref{Whitney}.\arabic*}]
\item[(1)]Given $Q\in \mathcal{W}$, there exists $x\in Q$ such that $Q'\in\mathcal{W}$ and $x\in Q'$ implies $Q'=Q$.
\item[(2)] There exist constants $A,B>1$ such that $$\diam(Q)\leq \frac{1}{A}\cdot \dist(Q,\partial U)\leq B\cdot \diam(Q).$$
\eenum
\edefn

As an easy consequence of the definitions we get
\blem\label{Whitney_lemma} Suppose the metric space $Y$ has a dyadic decomposition $\mathcal{D}$.
For any open set $U\subsetneq Y$, there exists an $A$-Whitney decomposition of $U$. The constant $B$ depends only on $\Gamma$ and $A$.
\elem
\bproof Let $A>10$.
For $\xi\in U$, choose $Q=Q_\xi\in\mathcal{D}$ such that $\xi\in Q$ is maximal (with respect to diameter) when $$\diam(Q)\leq \frac{1}{A}\cdot \di(\xi,\partial U).$$ This exists because of property (1) of \ref{dyadic}. Then by the maximality of $Q$, (3) and (5) of \ref{dyadic}, there exists $B_\Gamma$, such that $$\frac{1}{A}\cdot \di(\xi,\partial U)\leq B_\Gamma\cdot \diam(Q).$$ Then, $$\diam(Q)\leq \frac{1}{A-1}\cdot \di(Q,\partial U)\leq \frac{A\cdot B_\Gamma}{A-1}\cdot \diam(Q).$$This shows (2).

A subcollection of $\{Q_\xi\}_\xi$ may be chosen so that also (1) is satisfied. 
\eproof

\subsection{Hausdorff content and dimension} 
Let $Y$ be a metric space. Let $t\in [0,\infty)$. The $t$-Hausdorff content of a set $E\subset Y$ is defined as $$\mathcal{H}_\infty^t(E):=\inf\left\{\sum_i \diam(E_i)^t:E\subset \bigcup_i E_i\right\}.$$ The Hausdorff dimension of $E$ is defined as $$\text{dim}_{\mathcal{H}}(E):=\inf\,\{t\in [0,\infty):\mathcal{H}_\infty^t(E)=0\}.$$ The $t$-Hausdorff content is an outer measure. It is finite for bounded sets. It is not a Borel measure in general.

\subsection{CAT($-1$) spaces}\label{CAT_prelim} For the results on the spiraling of geodesics in \S \ref{sec:spiraling} we will have to assume curvature bounds for our hyperbolic spaces, more precisely, we consider manifolds $M$ of pinched negative sectional curvature, $-a^2\leq k \leq 1$. The methods however will not crucially depend on the smooth structure (see Remark \ref{general}). 

In \S \ref{sec:spiraling} we will use the Alexandrov `thin'-CAT($-1$) inequality for triangles in CAT($-1$) spaces and the `fat'-CBB($-a$) inequality for spaces with curvature bounded below (by $-a^2$), see for example \cite{BriHae} or \cite{AleKapPet}.
Another fact we will use is that for a CAT($-1$) space $X$, $$d_x(\xi,\eta)=e^{-(\xi|\eta)_x},$$ for any point $x\in X$ defines a visual metric on the visual metric (see \cite{Bou}), where $(\xi|\eta)_x$ is the Gromov product with base-point $x$.

\bdefn[Convex cocompact action]
We say that a group $\Gamma$ acts on a proper, geodesic, hyperbolic space $X$ convex cocompactly, if it acts properly and cocompactly on the convex hull (in $X$) of the limit set $\Lambda_\Gamma\subset\partial X$.
\edefn
\subsection{Dynamics of the $\Gamma$-action on the boundary}\label{CAT_prelim2} We will assume in \S \ref{sec:spiraling} that for a strictly negatively curved manifold $M$ the natural action of the fundamental group $\pi_1(M)$ on the visual boundary $\partial \tilde{M}$ (with the corresponding Patterson-Sullivan measure) of the Riemannian universal cover $\tilde{M}$ is ergodic.

\bdefn[Trail] Given $E\subset \tilde{M}$, define the \textit{trail}
$$T(x,E)=\{y\in\tilde{M}\mid \exists\,z\in\gamma_{xy},\;\text{such that}\,z\in E\},$$
\edefn
We will assume the following estimate for the distribution of orbits: there exists $n_0=n_0(M)>0$ such that \beqn\label{orbit_counting1}\#\{g\in\pi_1(M)\mid   gx \in B(x,n+n_0)\sm B(x,n)\}\sim e^{n\cdot v_{\pi_1(M)} },\eeqn and \begin{align}\label{orbit_counting2}\begin{split}\#&\{g \in\Gamma \mid gx\in T(x,B(y_0,L))\bigcap\left(B(x,n+n_0)\setminus B(x,n)\right)\}\\&\hspace{3.5cm}\approx \mu_{x}(S(x,B(y_0,L)))e^{n\cdot v_{\pi_1(M)}},\end{split}\end{align}  for $n$ large where $B(x,n)$ is a ball of radius $n$ in the universal cover $\tilde{M}$ of $M$, $L>0$ and $y_0\in\tilde{M}$ is such that $\gamma_{xy_0}(\infty)\in\Lambda_{\pi_1(M)}$ and $L$ and $\tilde{\rho}(x,y_0)$ are sufficiently smaller than $n$. The second inequality says that the density of orbit points in the intersection of an annulus with a cone is proportional to the measure of the shadow of the cone. These conditions are satisfied for example when the action of $\pi_1(M)$ on $M$ is convex cocompact or more generally, the quotient by $\Gamma$ of the `space of geodesic lines' admits a finite Bowen-Margulis-Sullivan measure for $\pi_1(M)$ (see Th\'{e}or\`{e}me 4.1.1 and Corollaire 2 in \cite{Rob}). 
Note that if $\pi_1(M)$ acts geometrically on $\tilde{M}$, the topological entropy of the geodesic flow in $M$, $v_{M}= v_{\pi_1(M)}$ (see \cite{Man}).

\section{Diophantine Approximation}\label{def:diophantine} 
In this section we describe a fine Diophantine approximation theory in a general setting, which we use later to study geodesic approximation in manifolds with negative curvature. A powerful technique used to obtain Hausdorff dimension lower bounds for `limsup type' sets which arise in Diophantine approximation are the regular systems of Baker and Schmidt \cite{BakSch}. This construction has been generalized to study cuspidal excursions of geodesics in hyperbolic manifolds by Meli\'{a}n and Pestana \cite{MP} in relation with Jarnik-Besicovitch theorems in the plane. One of our main innovations is to generalize the aforementioned work to allow for arbitrary open sets. 

\subsection{Generalized Jarnik-Besicovitch}
Let $X$ be a countable set and let $(A,\di,\mu)$ be a locally-compact metric space with a $D$-Ahlfors regular probability measure $\mu$ for $D>0$. Let $F:X \to \mathcal{P}(A)$ be a subset-valued map. Then the set \[
  E_F =\left\lbrace \xi\in A \;\middle|\;
  \begin{tabular}{@{}l@{}}
    $\exists\; \text{infinitely many}\;  x\in X$\; \text{such that} $\xi\in F(x)$
   \end{tabular}
  \right\rbrace
\] will be called $F$-\textit{approximable}. 
In the classical situation $X$ is a subset of a group acting on the space $A$ by homeomorphisms and the sets $F(x)$ are balls, for many applications however they will be more complicated sets (see Theorem \ref{version of 4.6} and \ref{spiral_01}).

We will now prove a general Jarnik-Besicovitch theorem. Let $$Dir:X\to\mathcal{P}(A)$$ be a (Dirichlet) function such that $x\in X$, $Dir(x)$ is a ball and,
\benum  
\item $\lim_{y} \text{rad}(Dir(y))\to 0$ (along every distinct sequence in $X$)
\item $A=E_{Dir}$ 
\item $\text{rad}(Dir(x))\leq 1$ for each $x\in X$. 
\eenum 
We will call the balls $Dir(x)$ Dirichlet balls.
Note the third property is not a real restriction, as the previous two properties always allow for modifying the set $X$ suitably, so that we get a system where the third property also holds. Nevertheless, we include it in the definition of $Dir$ for ease of exposition.

Let $$F:X\to \mathcal{P}(A)$$ be a (Jarnik-Besicovitch) function such that $F(x)\subset Dir(x)$ are open sets in $A$. The Hausdorff dimension of $E_F$ can be obtained in terms of the numbers introduced below which depend only on functions $Dir$, $F$ and the exponent of regularity of $\mu$, $D$.
For $x\in X$ write 
\[
\alpha_x:=\inf \left\lbrace\alpha\geq 1\;\middle|\; 
\begin{tabular}{@{}l@{}} $\text{there is a finite collection of mutually disjoint balls}\,\{B_i\}_i$\\$\text{such that for all}\; i,  \, \text{center}(B_i)\in F(x),$\\$ \text{rad}(B_i)=\text{rad}(Dir(x))^{\alpha},\;  \frac{1}{10^D}\cdot\mu(F(x))\leq\sum_i\mu(B_i)\leq \mu(F(x))$\end{tabular}\right\rbrace,
\] and 
$$\beta_x:=D\cdot\left(\alpha_x-\frac{\log(\mu(F(x)))}{\log(\mu(Dir(x)))}\right).$$ Since $\mu$ is a finite $D$-regular measure, $\alpha_x$ is well-defined. Indeed, to check this, one picks a compact subset of $F(x)$, close to it in measure, and applies the $5r$-covering theorem to a suitable collection of balls covering this compact set. We think of $\beta_x$ roughly, as a measure of how `far' $F(x)$ is from being a ball, indeed, it is zero when $F(x)$ is a ball. Set 
$$\overline{\alpha}_F:=\limsup_x \alpha_x,\,\underline{\beta}_F:=\liminf_x \beta_x,\,\underline{\alpha}_F=\liminf_x\alpha_x,\,\overline{\beta}_F=\limsup \beta_F.$$ 
Next for $k\in\R$, for each $\xi\in A$ consider the element $B_\xi$ in $Dir(X)$ of maximal radius containing $\xi$, such that the radius is bounded above by $1/k$. Call such a ball $k$-maximal. Let $B\subset A$ be a ball and consider any such cover $B=\bigcup_{\xi\in A} B_{\xi}$. For $0<c\leq 1$, $k\in\N$, and $B\subset A$, a ball, a collection of Dirichlet balls $\{B_{\xi_i}\}_i$ contained in $B$ is $(c,k)$-admissible if:
\begin{enumerate}
    \item $B_{\xi_i}=B(\xi_i,\text{rad}(B_{\xi_i}))$ are $ck$-maximal balls with $\xi_i\in B$, 
    \item $\bigcup_{\xi_i\in B} B_{\xi_i}\subset \frac{1}{2}B$,
    \item $1/ck\geq \text{rad}(B_{\xi_i})$,
    \item $c\mu(B)\leq \mu(\bigcup_i B_{\xi_i}).$
\end{enumerate}
Note that for $c$ small enough and $k$ large enough, a $(c,k)$-admissible collection of Dirichlet balls exists: one first obtains a finite cover of, say $\frac{1}{4}\overline{B}$ by Dirichlet balls and then applies the $5r$-covering theorem, to get a disjoint subcollection. Next, for a $(c,k)$-admissible collection $\{B_{\xi_i}\}_i$, write $\omega_{k}^{(i)}$ for the minimal element in $\N$, such that $1/ck\geq \text{rad}(B_{\xi_i})\geq c/k^{\omega^{(i)}_{k}}$, for each $\xi_i$. Note that $$\sup_{i}\omega^{(i)}_{k}<\infty,$$ since the collection $\{\xi_i\}_i$ is finite.

Write
\[
\omega_k(c,B):=\inf \left\lbrace\omega\geq 1\;\middle|\; 
\begin{tabular}{@{}l@{}}  $\omega=\sup_i \omega^{(i)}_{k}$, where $\{B_{\xi_i}\}_i$\\ is a $(c,k)$-admissible collection \end{tabular}\right\rbrace,
\]

and set $$\overline{\omega}(c,B)=\limsup_k \,\omega_k(c,B).$$ 
Set $$\overline{\omega}(c)=\sup\{\overline{\omega}(c,B)\mid B\,\text{is a ball in}\,A\},$$ and  $$\overline{\omega}=\inf_{c>0}\, \overline{\omega}(c).$$ Note that $\overline{\omega}$ only depends on $Dir$ and not on $F$. In our applications $\overline{\omega}$ will fortunately be rather easy to determine.

\bthm[Generalized Jarnik-Besicovitch]\label{density}Let $X$ be a countable set. Let $(A,\di,\mu)$ be a proper metric-measure space. Let a Dirichlet function $Dir$ and a Jarnik-Besicovitch function $F$ be given.
Assume 
$$\underline{\beta}_F>0,\; \overline{\alpha}_F<\infty$$ and there exists $0<c\leq 1$ such that $$\overline{\omega}(c)<\infty.$$
 Then, for $\alpha'>\overline{\alpha}_F$, $\beta'<\underline{\beta}_F$ and $\omega'>\overline{\omega}(c)$, $$\mathcal{H}^{d'}_{\infty}(E_F\cap B)\geq c(\alpha')\diam(B)^{d'},$$ where $d':=\frac{1}{\omega'}\cdot \left(\frac{\beta'+D}{\alpha'}\right)$ for any ball $B\subset A$. Moreover, $$\frac{1}{\overline{\omega}}\cdot\left(\frac{\underline{\beta}_F+D}{\overline{\alpha}_F}\right)\leq \;\dim_{\mathcal{H}}(E_F)\; \leq \;\frac{\overline{\beta}_F+D}{\underline{\alpha}_F}.$$
\ethm
\bproof
(\textbf{Lower bound.})\newline \textbf{Step 1}. Suppose that $B\subset A$ is a ball, and $k\geq k_B$, where $k_B$ is the smallest element in $\N$ such that $\omega_{k'}(c,B)<\omega'$, for all $k'\geq k_B$. Consider a collection of $ck$-maximal balls $\{B_{\xi_i}\}_i$ contained in $\frac{1}{2}B$, such that $c\mu(B)\leq \mu(\bigcup_i B_{\xi_i})$. Let $\hat{B}_{\xi_i}$ be the balls with same centre and radius $2/ck$. Applying the 5r-covering theorem to the collection $\{\hat{B}_{\xi_i}\}_i$, we get a disjoint subcollection again denoted $\{\hat{B}_{\xi_i}\}_i$ and write $$\tilde{\mathcal{F}}_{B,k}:=\{\text{Dir}(x_i)\mid \text{center}(\hat{B}_{\xi_i})=\text{center}(Dir(x_i)),\,1\leq i\leq i_{B,k}\}.$$ Note that $$\#\tilde{\mathcal{F}}_{B,k}\geq \frac{1}{C'}\cdot\mu(B)\cdot k^{D},$$ where $C'$ is a function of $D$.

Now we obtain approximations by suitable coverings to the sets in the collection $F(X)$. First note that there exists a collection of distinct points $\{\xi_i\}_i\subset F(x)$, $\di(\xi_i,\xi_j)\geq \text{rad}(Dir(x))^{\alpha_x}$ such that $$F(x)\subset\bigcup_i B(\xi_i,\text{rad}(Dir(x))^{\alpha_x}),$$ and upto constants depending on $D$ it holds,
$$\mu(F(x))\approx \sum_i \mu(B(\xi_i,\text{rad}(Dir(x))^{\alpha_x})).$$ Denote by $n_x$ the number of balls in the cover and let $\tilde{\beta}_x$ be such that $$\log(n_x)=\tilde{\beta}_x\log\left(\frac{1}{\text{rad}(Dir(x))}\right).$$
Then an elementary computation reveals 

$$\liminf_x \beta_x=\liminf_x \tilde{\beta}_x$$ Write $$\mathcal{I}_{x}= \{B(\xi_i,\text{rad}(Dir(x))^{\alpha_x})\}_i,$$ for the disjoint collection obtained above. We will call the collections $\{\mathcal{I}_x\}_{x\in X}$ $F$-colonies below. \newline Write 
\[
  \hat{\mathcal{F}}_{B,k}:=\left\lbrace B(\xi',r_{B,k}) \;\middle|\;
  \begin{tabular}{@{}l@{}}
     $B(\xi',\text{rad}(Dir(x))^{\alpha_{x}})\in \mathcal{I}_x$,   \\   $Dir(x)\in\tilde{\mathcal{F}}_{B,k}$,    \\     $r_{B,k}=\min\{\text{rad}(Dir(x))^{\alpha_x}/2,1/4ck\}$
   \end{tabular}
  \right\rbrace.
\] and 
$$\tilde{\beta}_k:= \inf\{\tilde{\beta}_x\mid \text{rad}(Dir(x))\leq 1/k\}.$$ Construct a subcollection $\mathcal{F}_{B,k}$ from $\hat{\mathcal{F}}_{B,k}$ by choosing $k^{\tilde{\beta}_k}$ balls in $\hat{\mathcal{F}}_{B,k}$, corresponding to $\mathcal{I}_x$ for each $Dir(x)\in \tilde{\mathcal{F}}_{B,k}$. Then the collection $\mathcal{F}_{B,k}$ satisfies the following properties. 
\benum
\item $\bigcup\{B'\mid B'\in \mathcal{F}_{B,k}\}\subset B,$
\item $\frac{1}{C}\cdot\mu(B)\cdot k^{D+\tilde{\beta}_k}\leq \#\mathcal{F}_{B,k}$, where $C>0$ is a constant depending only on $D$ and
\item If $B', B''\in \mathcal{F}_{B,k}$, $B'$ is a concentric to a ball $\tilde{B'}\in\mathcal{I}(x)$ and $B''$ is concentric to a ball $\tilde{B''}\in\mathcal{I}(y)$ with $x\neq y$, then $$\dist(B',B'')\geq 1/2ck,$$  
\eenum
where we abbreviate $\omega_{k}(c,B)$ as $\omega_k$.

\vspace{2.5mm}
\textbf{Step 2}. We will obtain a constant $0<\overline{C}<1$ depending only on $D$ such that if for any closed ball $B\subset A$ and any collection of open sets $\{U_i\}_i$, if $$\sum_i \diam(U_i)^{d'}\leq C\cdot\text{rad}(B)^{d'},$$ then $\{U_i\}_i$ can not cover $E_F$. A point in $E_F$ which does not lie in $\bigcup_i U_i$ will be obtained as the intersection of a nested sequence of closed sets constructed inductively. First for $k\in\N$, let us set $$\alpha_k=\sup\{\alpha_x\mid \text{rad}(Dir(x))\leq 1/k\}.$$

Now, start by writing $V_0:=B$, $\mathcal{J}_0=\{B\}$,  $k_0:=k_B$ and $l_0:=C^{\frac{1}{d'}}/k_0$. Suppose $V_0,V_1,\ldots, V_{i}$ have been defined and satisfy for $1\leq j\leq i$
\benum 
\item $l_{j}=k_{j}^{-\alpha'\omega_{k_j}}\cdot \mu(V_{j-1})^{-\frac{1}{d'}}$

\item $k_j$ satisfies $$k_{j}^{-\alpha'\omega_{k_j}}\cdot \mu(V_{j-1})^{-\frac{1}{d'}}< k_j^{-\alpha_{k_j}\omega_{k_j}}.$$ 
\item $V_j= \{B\mid B\in\mathcal{J}_j\}$ is a finite union of disjoint balls concentric to balls in $F$-colonies, and of radius $l_j$, such that $$\#\mathcal{J}_j\geq \frac{1}{c}\cdot k_j^{D+\tilde{\beta}_{k_j}}\cdot \mu(V_{j-1}),$$ where $c$ is a constant depending only on $v_\Gamma$.
\item For $B'\neq B''$ in the collection $\mathcal{J}_j$ above coming from distinct $F$-colonies, $$\dist(B',B'')\geq \frac{1}{2k_j}.$$
\item $V_j\cap U_r=\emptyset$ if $l_j<\diam(U_r)\leq l_{j-1}$.
\item $\mu(V_j)\geq c'\cdot  \l_j^{D}\cdot k_j^{D+\beta_{k_j}}\cdot\mu(V_{j-1})$ where $c'$ depends only on $D$.
\eenum
where $0<\overline{C}<1$ is the number to be specified below. 
Then $V_{i+1}$ is constructed in the following way. Let $k_{i+1}\geq k_i+1$, $k_{i+1}\geq k_{B'}$, for all $B'\in \mathcal{J}_{i}$, be a number large enough that $\alpha_{k_{i+1}}<\alpha'$, $\beta_{k_{i+1}}>\beta'$, $\omega_{k_{i+1}}<\omega'$, $k_{i+1}\geq 8/l_i$ and $$k_{i+1}^{-\alpha'\omega_{k_{i+1}}}\cdot \mu(V_{i})^{-\frac{1}{d'}}< k_{i+1}^{-\alpha_{k_{i+1}}\omega_{k_{i+1}}}.$$ Define $$l_{i+1}=\frac{1}{2}\cdot k_{i+1}^{-\alpha'\omega_{k_{i+1}}}\cdot \mu(V_{i})^{-\frac{1}{d'}}.$$
Given $B'\in \mathcal{J}_{i}$, use Step 1 to obtain the collection of balls $\mathcal{F}_{B',k_{i+1}}$. Let $\mathcal{J}'_{i+1}$ be collection of all such balls coming from each $B'\in\mathcal{J}_i$ with same centres but the smaller radius $2l_{i+1}$. Consider the collection $$\mathcal{U}_{i+1}=\{U_r\mid l_{i+1}<\diam(U_r)\leq l_i\}.$$ Then apply Lemma 2.3 of \cite{MP} as in \cite{MP} (which only uses Ahlfors-regularity) to the collection of open sets $\mathcal{U}_{i+1}$ and a collection $\mathcal{B}'_{i+1}$ of balls obtained by choosing exactly one ball corresponding to an $F$-colony contributing to the collection $\mathcal{J}'_{i+1}$ whenever such a ball intersects with a set of the collection $\mathcal{U}_{i+1}$ in a ball of radius greater than $l_{i+1}$. Note that distinct such balls are at least $(1/2k_{i+1})$-apart. Then we get the bound $$\#\mathcal{B}'_{i+1}\leq \overline{C}\cdot \tilde{C}\cdot k_{i+1}^{D}\cdot \left(\frac{1}{k_0}\right)^{d'}\cdot l_j^{D-d'}.$$ Then by induction hypothesis (6), $$\#\mathcal{B}'_{i+1}\leq \overline{C}\cdot \tilde{C}\cdot k_{i+1}^{D} \cdot l_j^{-d'} \frac{\mu(V_i)}{c'\cdot k_i^{\beta_{k_i}+D}\cdot \mu(V_{i-1})}.$$ Then, by induction hypothesis (1), 
$$\#\mathcal{B}'_{i+1}\leq \overline{C}\cdot \tilde{C}\cdot k_{i+1}^{D}\cdot \frac{\mu(V_i)}{c'\cdot k_i^{\beta_{k_i}-\beta'}}\leq \overline{C}\cdot \hat{C} \cdot k_{i+1}^{D}\cdot \mu(V_i),$$ where $\hat{C}$ depends only on $D$. Let $\mathcal{B}_{i+1}$ be the collection of all balls in $\mathcal{J}_{i+1}$ which intersects a set in $\mathcal{U}_{i+1}$ in a ball of radius greater than $l_{i+1}$. Since each Dir-colony contributes at most $k_{i+1}^{\beta_{k_{i+1}}}$ balls, we have that $$\#\mathcal{B}_{i+1}\leq \overline{C}\cdot \hat{C} \cdot k_{i+1}^{D+\beta_{k_{i+1}}}\cdot \mu(V_i).$$ Note that by properties of the collections $\mathcal{F}_{B,k}$ we also have that $$\#\mathcal{J}'_{i+1}\geq \frac{1}{C'}\cdot k_{i+1}^{D+\beta_{k_{i+1}}}\cdot \mu(V_i).$$ Now, $\overline{C}$ is chosen small enough depending only on $D$ such that $$\#(\mathcal{J}'_{i+1}\setminus \mathcal{B}_{i+1}) \geq \frac{1}{c}\cdot k_{i+1}^{D+\beta_{k_{i+1}}}\cdot \mu(V_i).$$ Define $\mathcal{J}_{i+1}$ to be the collection of balls which are concentric to the balls in $\mathcal{J}'_{i+1}\setminus \mathcal{B}_{i+1}$ and have radius $l_{i+1}$. Set $$V_{i+1}:=\bigcup\{B'\mid B'\in\mathcal{J}_{i+1}\}.$$ This finishes the induction step. This also shows that the intersection $\bigcap_i V_i\subset E_F$ is disjoint from the union $\bigcup_i U_i$ and thus the lower bound holds.

\vspace{2mm}
(\textbf{Upper bound.})\newline For $k\in\N$, write $$\alpha'_k=\inf \{\alpha_x\mid \text{rad}(Dir(x))\leq 1/2^k\},$$ $$\beta'_k=\sup \{\beta_x\mid \text{rad}(Dir(x))\leq 1/2^k\}.$$ For any $\alpha''<\underline{\alpha}_F$, $\beta''>\overline{\beta}_F$, and set $d'':=\left(\frac{\beta''+D}{\alpha''}\right)$. Choose $k_0$ large enough so that for all $k\geq k_0$, $\alpha'_k>\alpha''$ and $\beta'_k<\beta''$.
Now fix $k\geq k_0$. It suffices to show that given $0<\delta<1/2^{k_0}$, for all $d''>\frac{D+\beta''}{\alpha''}$, $$\mathcal{H}_{\delta}^{d'}(E_F)= O(\delta^{\,\alpha''d''-D-\beta''}).$$

Let $0<\delta<1/2^{k_0}$ be fixed now and let $i_0\in \N$ be smallest such that $1/2^{i_0}\leq \delta$. Next, given $\xi\in E_F$, choose $x_\xi^{(i_0)}$ such that $Dir(x_\xi^{(i_0)})$
has maximal diameter constrained to the conditions $\xi\in F(x_\xi^{(i_0)})$ and $\text{rad}(Dir(x_\xi^{(i_0)}))\leq 1/2^{i_0}$.

Then define 

$$S_i=\left\{\xi\in E_F\,\mid\, \frac{1}{2^{i+1}}<\text{rad}(Dir(x_\xi^{(i_0)}))\leq \frac{1}{2^{i}}\right\}.$$ 
Note that $$E_F=\bigcup_{i=i_0}^{\infty} S_i.$$

Apply the $5r$-covering theorem to the collection $\{Dir(x_\xi^{(i_0)})\}_{\xi\in S_i}$, to obtain a countable disjoint subcollection of balls $\mathcal{S}_i$ such that $$S_i\subset \bigcup_{\xi\in S_i} F(x_\xi^{(i_0)})\subset \bigcup_{\tilde{B}\in \mathcal{S}_i} 5\tilde{B}.$$ 
Apply again the covering theorem to $$\bigcup_{\xi\in S_i} F(x_\xi^{(i_0)})\subset \bigcup_{\xi\in S_i}\;\;\bigcup_{B\;\in\; \mathcal{I}_{x_\xi^{(i_0)}}} B,$$ where $\mathcal{I}_{x_\xi^{(i_0)}}$ are defined in the proof of the dimension lower bound above, to obtain an atmost countable collection of points $\xi_l\in S_i$ such that $$S_i\subset \bigcup_{\xi\in S_i} F(x_\xi^{(i_0)})\subset \bigcup_{\xi_l\in S_i}\;\;\bigcup_{B\;\in\; \mathcal{I}_{x_{\xi_l}^{(i_0)}}} 5B.$$ Note further that a volume comparison argument provides that for each $\xi_l$, the set $$\mathcal{N}_l=\{\tilde{B}\in\mathcal{S}_i\,\mid\, 5\tilde{B}\;\;\bigcap\;\; (\;\bigcup_{B\;\in\; \mathcal{I}_{x_{\xi_l}^{(i_0)}}} 5B\,)\neq \emptyset\}$$ satisfies $$1\leq \#\mathcal{N}_l\leq C,$$ where $C$ is an absolute constant depending only on regularity of $\mu$.

Then 

\beqnn \begin{split}\mathcal{H}_{\frac{1}{k}}^{d''}(S_i)& \leq \sum_l \sum_{B\;\in\; \mathcal{I}_{x_{\xi_l}^{(i_0)}}} \frac{1}{2^{i\alpha'_id''}}  \\& \leq \sum_l 2^{i\beta'_i}\frac{1}{2^{i\alpha'_id''}} \\& \leq \sum_l 2^{i(\beta''-\alpha''d'')}\quad (\text{since} \; i>k_0)\\ & \leq \frac{1}{C}\cdot \sum_{\mathcal{S}_i} 2^{i(\beta''-\alpha''d'')}\quad (\text{by the estimate on}\; \#\mathcal{N}_l)\\& \leq \frac{C'}{C}\cdot 2^{i\cdot(\beta''-\alpha''d''+D)},\end{split}\eeqnn where the last inequality followed from $\#\mathcal{S}_i \leq C'\cdot 2^{i\cdot D}$, which uses the regularity of $\mu$ and disjointedness of the elements of $\mathcal{S}_i$. Thus, $$\mathcal{H}_{\frac{1}{k}}^{d''}(E_F)\leq \sum_{i=i_0}^{\infty} 2^{i\cdot(\beta''-\alpha''d''+D)}.$$

The claim follows.
\eproof

\brem We have the following remarks on the previous theorem.
\benum 
\item For a similar statement in $\mathbb{R}^n$ with the Lebesgue measure, giving a Hausdorff lower bound in terms of Riesz energies, see Persson \cite{Persson}.
\item In case $A=E_{Dir}$ the existence of admissible $(c,k)$ covers follows by definition, for some $0<c\leq 1$ and $k\in\N$, as noted above. However, note that $\mu(A\setminus E_{Dir})=0$ with the assumption $\overline{\omega}(c)<\infty$ for some $0<c\leq 1$ also provides such covers and is sufficient for the proof to work. 
\item One can modify the above method to generalise to the case where the Dirichlet cover is given in terms of more general open sets with some regularity condition to get a slightly more general Jarnik-Besicovitch statement. As mentioned before, this circle of problems has attracted attention recently. 

\item In our applications we always have $\overline{\omega}=1$. 

\eenum
\erem

\subsection{Large intersections}
Consider a compact metric space $A$ as above with a dyadic decomposition $\mathcal{D}$.

Let $s>0$ be given. We define first the class 
\[
  \mathcal{G}_s:=\left\lbrace E\subset A \;\middle|\;
  \begin{tabular}{@{}l@{}}
     $E\;\text{is}\;G_\delta\;\text{and}\;\mathcal{H}_{\infty}^{s}(F\cap Q)\geq c\cdot \diam(Q)^s$\\ for some  $c>0$, for all $Q\in\mathcal{D}$ 
   \end{tabular}
  \right\rbrace.
\]
The limsup sets we will study in this paper will belong to such a collection $\mathcal{G}_s$ for suitable $s$.
In this section we prove a version of Falconer's large intersection property from \cite{Fal} for $\mathcal{G}_s$.

\bdefn[Parent]\label{parent}
Let $F\subset A$. By a parent of $F$, we mean a set $Q\in\mathcal{D}$ which is a minimal (with respect to diameter) for sets of $\mathcal{D}$ containing $F$. We denote it by $\hat{F}$. 
\edefn

\bthm\label{intersect_main}
If $E\in\mathcal{G}_s$ for some $0<s<D$, then we have for all $t\in[0,s)$ and $U\subset A$ open that, \beqn\label{main_eqn}\mathcal{H}_{\infty}^{t}(E\cap U)\geq \frac{\mathcal{H}_{\infty}^s(U)}{\diam(U)^{s-t}}.\eeqn
\ethm

\bproof
Note that to estimate the left hand side of inequality \eqref{main_eqn}, it suffices to consider coverings by closed sets.
Let $\mathcal{I}=\{E_i\}_i$ be a covering of $F\cap U$ by closed sets. We wish to bound from below the sum $\underset{\mathcal{I}}\sum\diam(E_i)^t$. We may assume without loss of generality that $\diam(E_i)\leq \diam(U)$.

Set $r>0$, be such that $$r=\alpha\cdot \diam(U),$$ and 
$$\mathcal{I}_{\alpha}:=\{E_i\in\mathcal{I}:\diam(E_i)\geq r\},$$ where $0<\alpha=\alpha(Q,s,t)<1$ will be specified later. 

Let $\mathcal{W}_{\alpha}'$ be a Whitney decomposition (see Definition \ref{Whitney}) of the open set $$U_{\alpha}:= U\setminus  \bigcup_{\mathcal{I}_{\alpha}}E_i,$$ for some $A>10$ (note that this exists in the case $\mathcal{I}_\alpha=\emptyset$ because $U\neq\partial X$).

For every $Q\in\mathcal{W}_{\alpha}'$, for which there exists $E'\in\mathcal{I}\setminus\mathcal{I}_{\alpha}$, such that a parent $\hat{E'}$ of $B(E',\frac{1}{2}|E'|)$ contains $Q$, we replace $Q$ by $\hat{E'}$ in the collection $\mathcal{W}_{\alpha}'$, and call the new collection $\mathcal{W}_{\alpha}$.

Note that $\mathcal{I}_{\alpha}\cup \mathcal{W}_{\alpha}$ is a cover for $U$. Also note that if $Q\in\mathcal{W}_{\alpha}\bigcap \mathcal{W}_{\alpha}'$, then for $E'\in\mathcal{I}\setminus\mathcal{I}_{\alpha}$ such that $\hat{E'}\bigcap Q\neq \emptyset$, we have $$\diam(Q)\geq \frac{1}{2}\diam(E').$$ Moreover, if $Q,Q'\in\mathcal{W}_{\alpha}\bigcap \mathcal{W}_{\alpha}'$ are such that $Q'\bigcap E'\neq\emptyset\neq Q\bigcap E'$ for some $E'\in\mathcal{I}\setminus\mathcal{I}_{\alpha}$, then 
$$\di(Q',\partial U_{\alpha})\leq \di(Q,Q')+\di(Q,\partial U).$$ Then by Definition \ref{Whitney} and Lemma \ref{dyadic_overlap}, there exists a constant $M$ such that
$$\#\{Q\in \mathcal{W}_{\alpha}\bigcap \mathcal{W}_{\alpha}':Q\bigcap E'\neq\emptyset\}\leq M.$$

\beqnn
\begin{split}
\sum_{\mathcal{I}}\diam(E_i)^t &
=\sum_{\mathcal{I}_{\alpha}}\diam(E_i)^t + \sum_{\mathcal{I}\setminus\mathcal{I}_{\alpha}}\diam(E_i)^t \\ &
\geq \diam(U)^{t-s}\sum_{\mathcal{I}_{\alpha}}\diam(E_i)^s  +(\alpha\diam(U))^{t-s} \sum_{\mathcal{I}\setminus\mathcal{I}_{\alpha}}\diam(E_i)^s \\ &
\geq \diam(U)^{t-s}\sum_{\mathcal{I}_{\alpha}}\diam(E_i)^s \\& +\frac{1}{2} \cdot \left(\frac{1}{c^s}(\alpha\cdot\diam(U))^{t-s}\sum_{\mathcal{W}_{\alpha}\setminus\mathcal{W}_{\alpha}'}\diam(\hat{E_i})^s\right)  \\&  + \frac{1}{2}\cdot \left(\frac{1}{M}(\alpha\cdot\diam(U))^{t-s}\sum_{Q\in\mathcal{W}_{\alpha}\cap\mathcal{W}_{\alpha}'}\sum_{\substack{E_i\in \mathcal{I}\setminus\mathcal{I}_{\alpha} \\E_i\cap Q\neq \emptyset}}\diam(E_i)^s\right) \\&
\geq \diam(U)^{t-s}\sum_{\mathcal{I}_{\alpha}}\diam(E_i)^s  + \frac{1}{2c'}(\alpha\cdot\diam(U))^{t-s}\sum_{Q\in\mathcal{W}_{\alpha}} \diam(Q)^s\\&
\geq \diam(U)^{t-s}\sum_{\mathcal{I}_{\alpha}\cup \mathcal{W}_{\alpha}}\diam(Q)^s\\&
\geq \frac{\mathcal{H}_{\infty}^s(U)}{\diam(U)^{s-t}}
,\end{split}
\eeqnn
for $\alpha>0$ small enough (which can be suitably chosen since $t<s$).

\eproof

\blem\label{intersect_2}
For $Q\in\mathcal{D}$, $\mathcal{H}_{\infty}^s(Q)\approx \diam(Q)^s$, for $s\in [0,D)$. The comparability constants depend only on $s$.
\elem

\bproof  Let $\{E_i\}_i$ be a covering of $Q$.
Let $\{Q_i\}_i\subset \mathcal{D}$ be the collection of respective parents. Then, \beqn\begin{split}\sum_i \diam(E_i)^s &\geq \frac{1}{c^s}\cdot\sum_{i}\diam(Q_i)^s =\frac{1}{c^s}\cdot\sum_{i}\diam(Q_i)^{D}\cdot \diam(Q_i)^{s-D}\\&\geq \frac{1}{c'^s}\cdot \sum_{i}\mu(Q_i)\cdot \diam(Q)^{s-D}\geq \frac{1}{c'^s}\cdot \mu(Q)\cdot\diam(Q)^{s-D} \\& \geq \frac{1}{c''^s}\cdot\diam(Q)^s, \end{split}\eeqn gives the claim. 
\eproof
\brem
Note that Lemma \ref{intersect_main}
 and \ref{intersect_2} imply that $\mathcal{G}_s\subset \mathcal{G}_t$ for $0\leq t<s$.\erem

\bdefn
We define for each $t\geq 0$, $$\hat{\mathcal{H}}_{\infty}^t(F):= \sup_{s \downarrow t}\mathcal{H}_{\infty}^s(F),$$ for every $F\subset A$.
\edefn

Below we define a version of the increasing sets property suitable to our situation. For details see for example Rogers \cite{Rog} or Howroyd \cite{How}.
\bdefn[Increasing sets property]
An outer measure $m$ satisfies the increasing sets property if for any collection $\{F_i\}_i$ of nested increasing sets we have that $$m\left(\bigcup_i F_i\right)=\lim_i m(F_i).$$
\edefn

\brem
The set function $\hat{\mathcal{H}}_{\infty}^t$ is an outer measure. It is bounded from above by $\mathcal{H}_{\infty}^t$ (up to a constant depending on the diameter of $A$; if $\diam(A)<1$, the constant can be taken to be one). It satisfies the aforementioned version of the increasing sets property; see for example Howroyd \cite{How} (page 29, Corollary 8.2). 
\erem

\blem\label{main}
Given $0<s<D$ and $E\in\mathcal{G}_s$, for all $t\in [0,s)$ we have $$\hat{\mathcal{H}}_{\infty}^t(E\cap U)= \hat{\mathcal{H}}_{\infty}^t(U).$$
\elem

\bproof
It is clear that $\hat{\mathcal{H}}_{\infty}^t(E\cap U)\leq \hat{\mathcal{H}}_{\infty}^t(U).$ For the other inequality we note that for $s>t$, and any $s>s'>t$, $$\hat{\mathcal{H}}_{\infty}^t(E\cap U)\geq \mathcal{H}_{\infty}^{s'}(E\cap U)\geq \frac{\mathcal{H}_{\infty}^s(U)}{\diam(U)^{s-s'}},$$ where the third inequality is Lemma \ref{intersect_main}. The claim follows by taking the limit of supremums as $s\downarrow t$.
\eproof

\bdefn[Metrically dense]Let $t\geq 0.$ A set $F\subset A$ is called $t$-metrically dense if for each open set $U\subset \partial X$ open, $$\hat{\mathcal{H}}_{\infty}^t(F\cap U)= \hat{\mathcal{H}}_{\infty}^t(U).$$
\edefn

We are now ready for 
\blem\label{compact} For all $t\in(0,D)$ the following holds. Let $\{F_i\}_{i\in\N}$ be a collection of $t$-metrically dense $G_\delta$. Let $U$ be an open set. Then
$$\mathcal{H}_{\infty}^t\left(\bigcap_{i}F_i\cap U\right)\gtrsim\hat{\mathcal{H}}_{\infty}^t(U).$$ It follows that when for $Q\in\mathcal{D}$, $$\mathcal{H}_{\infty}^t\left(\bigcap_{i}F_i\cap Q\right)\gtrsim\mathcal{H}_{\infty}^t(Q),$$ and $$\hat{\mathcal{H}}_{\infty}^t\left(\bigcap_{i}F_i\cap U\right)=\hat{\mathcal{H}}_{\infty}^t(U).$$ \elem

\bproof
The proof follows Falconer's argument in Lemma 4 of \cite{Fal}. First assume that $F_i$ are a sequence of decreasing open sets. Fix $\epsilon>0$ small.  Set $U_0:=U$. Then there exists a collection of numbers $\epsilon_i$ such that $$U_i:=\{x\in F_i\cap U_{i-1}:\di(x,\partial (F_i\cap U_{i-1}))> \epsilon_i\},$$ for $i\geq 1$ satisfy $$\hat{\mathcal{H}}_{\infty}^t(U_i)>\hat{\mathcal{H}}_{\infty}^t(U)-\epsilon,$$ for all $i$.

Observe that $\overline{U}_i\subset F_i\cap U$. Let $\{E_j\}_j$ be a covering of $\bigcap_i \overline{U}_i$. Then $\{\Int(\hat{E_j})\}_j$ is an open covering of $\bigcap_i \overline{U}_i$. There exists $k$ such that $\overline{U}_k\subset \bigcup_j \Int(\hat{E_j})$. Thus 

$$\sum_i \diam(E_j)^t\geq \sum_i \diam(E_j)^s\gtrsim \mathcal{H}_{\infty}^s(U_k)>\hat{\mathcal{H}}_{\infty}^t(U)-2\epsilon,$$ for $s>t$ small enough. The claim follows by taking infimum over all such coverings of $\bigcap_i \overline{U}_i$ and letting $\epsilon$ go to zero.

The second inequality follows from Lemma \ref{intersect_2} and the continuity of the exponential function. The third equality follows from Lemma \ref{main} and the previous inequalities in the claim. The general case follows as argued in Lemma 4 of \cite{Fal}.
\eproof

We have thus proved the following theorem.
\bthm[Large intersection property]\label{LIP} Let $(A,\di,\mu)$ be a proper, $D$-Ahlfors regular metric space, with a dyadic decomposition $\mathcal{D}$. Let $0<s<D$. Let $\mathcal{F}^s$ be a collection of $G_\delta$ sets. The following are equivalent.
\benum
\item For each $E\in\mathcal{F}^s$, and each $0\leq t <s$ and $Q\in\mathcal{D}$, it holds $$\mathcal{H}_\infty^t(E\cap Q)\geq c_{t,E}\mathcal{H}_\infty^t(Q),$$ for some $c_{t,E}>0$.

\item For each countable collection $\{E_i\}_i\subset\mathcal{F}^s$, and each $0<t<s$ and $U\subset\partial X$ open such that $U$ has a Whitney decomposition, $$\hat{\mathcal{H}}_\infty^t\left(\bigcap_i E_i\cap U\right)=\hat{\mathcal{H}}_\infty^t(U).$$ 
\eenum
\ethm
In the case $A$ is compact, the result follows from Lemma \ref{compact}. The result holds when $A$ is proper as $A$ can be exhausted by compact subsets with corresponding dyadic decompositions.
\brem We note that:
\benum
\item It is clear that if $E\in\mathcal{G}^s$, then $\text{dim}_{\mathcal{H}}(E)\geq s$.
\item It also holds that for each $E\in\mathcal{G}^s$, each $0\leq t<s$ and $U\subset A$ open, $$\mathcal{H}_\infty^t(E\cap U)\geq c'_{t,E}\mathcal{H}_\infty^t(U),$$ for some $c'_{t,E}>0$. To see this one compares a given covering of $E\cap U$ with the Whitney decomposition of $U$ as in Lemma \ref{intersect_main}. We do not present details as we do not require this later.
\eenum
\erem

\subsection{A Borel-Cantelli lemma}
We have the following version of the classical Borel-Cantelli lemma  tailored for our applications.
\blem[Borel-Cantelli]\label{Borel_Cant} Let $X$ be a countable set. Let $(A,\di,\mu)$ be a proper metric space, with a $D$-regular probability measure $\mu$.
Let $Dir:X\to \mathcal{P}(A)$ be a function, such that there exist $a>1$ and $b\geq 1$ for which \benum
\item $Dir(x)$ is a ball of radius bounded by $a$. \vspace{0.2cm}\item  $\frac{1}{b}\cdot e^{n\cdot D}\leq  \#\{x\in X \mid \frac{1}{a}\cdot e^{-n}\leq \text{rad}(Dir(x))\leq a\cdot e^{-n}\}\leq b \cdot e^{n\cdot D},$ and \vspace{0.2cm}\item given a point $\xi \in A$, at most $b$ of the balls corresponding to $\{x\in X \mid \frac{1}{a}\cdot e^{-n}\leq \text{rad}(Dir(x))\leq a\cdot e^{-n}\}$ contain it. \vspace{0.2cm}\item $\frac{1}{b}\mu(X)\leq \sum \{\mu(Dir(x)) \mid x\in X,\, \frac{1}{a}\cdot e^{-n}\leq \text{rad}(Dir(x))\leq a\cdot e^{-n}\}\leq b\mu(X).$ \eenum \vspace{0.2cm}
Let $F:X\to \mathcal{P}(A)$ be a function such that $F(x)\subset Dir(x)$ and with the notation $$i_x:=-\log(\text{rad}(Dir(x)))),$$ we have that
\benum \item $F(x)=\bigcup_{B\in\mathcal{I}_x} B,$ is a union of balls in $A$, where \vspace{0.2cm}\item $\frac{1}{c}\cdot e^{n(i_x)}\leq \#\mathcal{I}_x\leq c \cdot e^{n(i_x)},$ and \vspace{0.2cm}\item $\frac{1}{c} \cdot e^{-(i_x+\varphi(i_x))}\leq \text{rad}(B)\leq c \cdot e^{-(i_x+\varphi(i_x))},$ for all $B\in\mathcal{I}_x$ 
\eenum \vspace{0.2cm} where $c\geq 1$ and $n,\varphi:[0,\infty)\to \R$, are positive \textit{Lipschitz} functions.
The set $E_F$ then has measure zero if $$\int_{1}^\infty e^{-(\varphi(t)D-n(t))}dt$$ converges. Moreover, if given $\xi \in A$ and $x\in X$, at most $c$ balls of $F(x)$ contain $x$, then the set $E_F$ has positive measure if $$\int_{1}^\infty e^{-(\varphi(t)D-n(t))}dt$$ diverges.
\elem
\bproof Write $$A_k :=\bigcup_{x: i_x\in (k-\log(a),k+\log(a))}\bigcup_{B\in \mathcal{I}_x} B.$$ Then $$E_{F}= \lim_i\sup A_i.$$ It is clear since $f$ is Lipschitz that convergence of the integral implies that $E_F$ has zero measure. 

For the second claim, we show that a suitable subsequence of $A_i$'s are quasi-independent. Let $i,j\in\N$ be such that $i+f(i)<j$. Then given a ball $B$ lying in $A_i$, the number of balls making $A_j$ intersecting it is bounded above by (a constant multiple of) $$e^{n(j)}\cdot e^{-(i+\varphi(i))\cdot D}\cdot e^{j\cdot D}.$$ Thus $$\mu(A_i\cap A_j)\leq e^{n(j)}\cdot e^{-(i+\varphi(i))\cdot D}e^{j\cdot D}\cdot e^{-(j+\varphi(j))\cdot D}\cdot \#\mathcal{A}_i\lesssim \mu(A_i)\cdot \mu(A_j),$$ and the unstated constants in above inequalities depend only on $a$, $b$, $c$ and $c'$. It follows that $E_{F}$ has positive measure when the integral diverges. The convergence case is a standard application of the Borel-Cantelli lemma and the construction.
\eproof

\subsection{Summary}
We summarise the dimension results of this section in the following theorem.

\bthm 
Let $(A,\di,\mu)$ be a proper metric-measure space with the $D$-Ahlfors regular measure $\mu$, and $X$ be a countable set. 

Assume given a function $Dir:X\to \mathcal{P}(A)$, where $Dir(x)$ are balls with $\text{rad}(Dir(x))$ accumulating at zero, and the Dirichlet type statement $$A=\limsup_{x\in X} Dir(x).$$ 
\benum
\item Given any other function $F:X\to \mathcal{P}(A)$ such that $F(x)\subset Dir(x)$ are open for all $x\in X$, a Hausdorff-dimension lower bound $d_F$ for the Jarnik-Besicovitch set $$E_F=\limsup_{x\in X} F(x)$$ can be obtained (in terms of the asymptotic behaviour of F), as well as the density condition $$\hat{\mathcal{H}}^d_{\infty}(E_F\cap U)=\hat{\mathcal{H}}^d_{\infty}(U)$$ for any $0\leq d<d_F$ and open set $U\subset A$, whenever $d_F>0$, where $\hat{\mathcal{H}}^d_{\infty}$ is a suitable modification of the $d$-Hausdorff content.

\item If $\mathcal{F}=\{F_i\}_i$ is a (countable) sequence of functions such that $F_i(x)\subset Dir(x)$ are open for all $x\in X$, then $$\dim_{\mathcal{H}}\,\left(\bigcap_i E_{F_i}\right)=\inf d_{F_i}:=d_{\mathcal{F}}.$$ The \textit{Liouville} set $E_{\mathcal{F}}:=\bigcap_i E_{F_i}$ also satisfies $$\hat{\mathcal{H}}^d_{\infty}(E_{\mathcal{F}}\cap U)=\hat{\mathcal{H}}^d_{\infty}(U),$$ for any $0\leq d < d_{\mathcal{F}}$ if $d_{\mathcal{F}}>0$. If $d_{\mathcal{F}}=0$, but $d_{F_i}>0$ for all $i$, the same holds for $d=0$ (by the Baire-category theorem).
\eenum
\ethm

\section{Shrinking targets and spiral traps}\label{sec:spiraling}
This section is devoted to establishing fine logarithm law type results for geodesics in negative curvature, especially Hausdorff dimension results for `spiraling' phenomena.

\bdefn\label{main_properties} Let $X$ be a CAT($-1$) space and $\overline{X}$ its compactification with Gromov's boundary. Let $\Gamma$ be a discrete group of isometries of $X$ acting properly on $X$. We will say the $\Gamma$ action on $X$ has property (i), $i\in\{1,2,3,4\}$, if it satisfies respectively
\benum

\item The action of $\Gamma$ on $\partial X$ equipped with the corresponding Patterson-Sullivan density is ergodic.
\item The orbit counting estimates \eqref{orbit_counting1} and \eqref{orbit_counting2} from \S \ref{CAT_prelim2} hold.
\item The convex hull of the limit set $\Lambda_\Gamma$ in $X$ admits a cocompact, proper group action by a discrete subgroup of isometries (not necessarily $\Gamma$).
\item $\Lambda_{\Gamma}=\partial X$, where $\Lambda_\Gamma$ is the limit set of an orbit of $\Gamma$ in $\overline{X}$.
\eenum
\edefn
By the convex hull of a set $E\subset \partial X$, we mean the convex hull in $X$ of the union of all geodesic lines with their end points in $E$.

\subsection{The unit tangent space at a point and the visual boundary}We denote by $SM$ the unit tangent bundle of $M$ and by $SM_{x_0}$ the fiber over $x_0$. The unit tangent sphere $SM_{x_0}$ is equipped with the visual metric $d_{\tilde{x}_0}$ (see \S \ref{CAT_prelim}) by its identification to the visual boundary and with the Patterson-Sullivan measure corresponding to $\pi_1(M)$, for a chosen lift $\tilde{x}_0$ of $x_0$ in $\tilde{M}$. This metric-measure structure is better suited for studying the asymptotic properties of geodesics in variable curvature.

\subsection{Zero-one laws}

We now consider spiraling of long geodesic pieces into fixed neighbourhoods of a fixed closed, totally geodesic submanifold, a phenomenon we refer to as a \emph{spiral trap}. This theorem, namely Theorem \ref{version of 4.6} is due to  Hersonsky and Paulin (Theorem 4.6 in \cite{HP3}). We provide a different argument where we apply the results of \S \ref{def:diophantine}. We will also use part of the argument later for Hausdorff dimension computations.

\bthm\label{version of 4.6}
Let $\epsilon>0$ be fixed.
Let $M$ be a manifold of pinched negative sectional curvature $-a^2\leq k\leq -1$, $a\geq 1$, such that the covering action of $\Gamma:=\pi_1(M)$ on the convex hull of $\Lambda_{\Gamma}$ in $\tilde{M}$ has properties (1) and (2). Let $N$ be a compact, convex submanifold such that $\Gamma_N:=\pi_1(N)\xhookrightarrow{}{} \Gamma$ is injective, and $\Gamma_N$ is a non-trivial subgroup of $\Gamma$, such that $v_{\Gamma_N}<v_{\Gamma}$. Let $f$ be a positive Lipschitz function. Let $x_0\in M$. Then the set 
\[
  E^{f}_{N} =\left\lbrace v\in SM_{x_0} \;\middle|\;
  \begin{tabular}{@{}l@{}}
    $\exists\; \text{positive times}\; t_n\to\infty\;\text{such that}$\\ $\gamma_v(t_n,t_n+f(t_n))\subset B(N,\epsilon)$
   \end{tabular} 
  \right\rbrace
\] has full (resp. zero) measure if the integral $$\int_{1}^\infty e^{-f(t)(v_\Gamma-v_{\Gamma_N})}dt$$ diverges (resp. converges), where $\Gamma:=\pi_1(M)$ and $\gamma_v$ is the geodesic at $x_0$ at time zero with direction $v$.  
\ethm

\bproof
\textbf{Convergence.} Let $F_M$ be a fundamental domain for the action of $\Gamma:=\pi_1(M)$ on $\tilde{M}$. Let $\tilde{x}_0$ and $\tilde{N}_0$ be the components of preimages of $x_0$ and $N$ intersecting $F_M$ non-trivially. In this proof and below we will use the abbreviations $$\rho_g:=\tilde{\rho}(\tilde{x}_0,g\tilde{x}_0)$$ and $$\tilde{N}_g:=g\tilde{N}_0,$$ where $g\in\Gamma$.

First let $\tilde{\gamma}_v$ be the lift of the geodesic corresponding to a direction $v\in E^f_N$ starting at $\tilde{x}_0$. Let $t_n\to\infty$ be a sequence of times and $g_n'\in\Gamma$ be a sequence of isometries such that $$\tilde{\gamma}_{v}(t_n,t_n+f(t_n))\subset B(g_n'\tilde{N}_0,\epsilon).$$ Let $z_0$ be the nearest point projection from $\tilde{x}_0$ to the fundamental domain $F_N\subset \tilde{N}_0$ of $N$ in $F_M$. Let $g_n'$ be in the coset $[g_n]$ of $\Gamma/\Gamma_{N}$ (where $\text{Stab}(\partial\tilde{N}_0)=:\Gamma_N\simeq \pi_1(N)$) so that $g_n'=g_nh_n$ for some $h_n\in\Gamma_N$, where $g_n$ is such that 
$$
\tilde{\rho}(\tilde{x}_0,g_nz_0)=\dist_{\tilde{\rho}}(\tilde{x}_0,g_n\Gamma_N z_0)
$$

Then we have that \beqn\label{21}\dist_{\tilde{\rho}}(\tilde{\gamma}_v(t_n+f(t_n)),g_n\Gamma_N(z_0))\lesssim \tilde{\rho}(\tilde{\gamma}_v(t_n+f(t_n)),g_nh_nz_0)\leq c\eeqn (where $c$ is a positive depending on data). Write $\xi$ for the end point of the geodesic ray $\tilde{\gamma}_v$.
Then there exists a point $\xi^{g_nh_n}\in \partial \tilde{N}_{g_n}\subset \partial \tilde{M}$ in the shadow $S_{\tilde{N}_{g_n}}(g_nz_0,B_{\tilde{N}_{g_n}}(g_nh_nz_0,R_N))$ for some $R>0$ depending on data such that (since the triangle $[\xi,\tilde{x}_0,\xi^{g_nh_n}]$ is $a$-fat) \beqn\label{22}d_{\tilde{x}_0}(\xi,\xi^{g_nh_n})\leq c_2e^{-(t_n+f(t_n))},\eeqn where $B_{\tilde{N}_{g_n}}(\cdot,\cdot)$ and
$S_{\tilde{N}_{g_n}}(\cdot,\cdot)$ are respectively used to denote balls in $\tilde{N}_{g_n}$ and shadows in the boundary of the embedded space $\tilde{N}_{g_n}\subset \tilde{M}$ of balls in $\tilde{N}_{g_n}$. The last inequality holds because $g_nz_0$ is $\delta_M$-close to $\gamma_{\xi^{g_nh_n}}$ (which follows from the convexity of $N$ and by hyperbolicity of $\tilde{M}$).

Let $z_n$ be a nearest point from $\tilde{\gamma}_v(t_n+f(t_n))$ to the geodesic $\gamma_{\xi,g_nz_0}$. Then the CAT($-1$) inequality applied to the triangle $[x_0,\xi,g_nz_0]$ gives \beqn\label{23}
\tilde{\rho}(\tilde{\gamma}_v(t_n+f(t_n)),z_n)\leq c_1,
\eeqn where $c_1$ depends on the data (which follows by noting from hyperbolicity that $\dist_{\tilde{\rho}}(g_nz_0,\gamma_{\xi^{g_nh_n}})$ is bounded above, by thin-ness of the triangle $[\xi,x_0,\xi^{g_nh_n}]$ and using the triangle inequality).

By \eqref{21} and \eqref{23}, we get $$\tilde{\rho}(z_n,g_nh_nz_0)\leq c_3$$ where $c_3$ is a positive depending on data and thus (since the triangle $[\xi,g_nz_0,\xi^{g_nh_n}]$ is $a$-fat) \beqn\label{24}d_{g_nz_0}(\xi,\xi^{g_nh_n})\leq c_4e^{-\rho_{h_n}},\eeqn where $d_{g_nz_0}$ is the visual metric on $\partial\tilde{M}$ from basepoint $g_nz_0$.

Let $h_n'\in\Gamma_N$ be such that $g_nh_n'z_0$ is a nearest orbit point for $\tilde{\gamma}_v(t_n)$. Then from \eqref{24}, 
\beqn d_{g_nz_0}(\xi,\xi^{g_nh_n})\leq c_4e^{-(\rho_{h_n'}+(\rho_{h_n}-\rho_{h_n'}))} \leq c_5e^{-(\rho_{h_n'}+f(\rho_{g_nh_n'}))},\eeqn where $c_5$ is a constant depending only on data. 

Therefore we have that \[ 
E^f_N \subset \limsup \left\lbrace B_{gz_0}(\eta,c_5e^{-(\rho_{h}+f(\rho_{gh}))})\;\middle|\;
\begin{tabular}{@{}l@{}} $[g]\in\Gamma/\Gamma_N,h\in \Gamma_N,$\\ $\eta\in S(gz_0,B_N(ghz_0,R_N))\bigcap\partial\tilde{N}_g$ \end{tabular} \right\rbrace
\] where $B_{gz_0}(\cdot,\cdot)$ is used to denote a ball in the embedded space $\partial\tilde{N}_{g}\subset\partial\tilde{M}$ with visual metric $d_{gz_0}$. 

Applying the $5r$-covering theorem to the collection $$\left\{B_{gz_0}(\eta,c_5e^{-(\rho_{h}+f(\rho_{gh}))})\mid\eta\in S(gz_g,B_N(ghz_0,R_N))\bigcap\partial\tilde{N}_g\right\},$$ we get a finite collection $\{B_{gz_0}(\eta_i,c_5e^{-(\rho_{h}+f(\rho_{gh}))})\}$ of disjoint balls such that concentric balls with five times the radius of the balls in the collection, cover the original collection. By Ahlfors regularity (in the metric space $\partial\tilde{N}_g$), we know that the number of balls in this subcollection lies between constant positive multiples of $e^{f(\rho_{gh})v_{\Gamma_N}}$. Call the set of centres $\mathcal{J}_{gh}$.
Note then that a finite number of balls centred at $\eta_i$ of radius $c_7e^{-(\rho_{gh}+f(\rho_{gh}))}$ (in the metric $d_{\tilde{x}_0}$) for large enough $c_7$ depending only on data, with varying $g$ and $h$ (using a compactness argument giving an upper bound on the number of lifts $\tilde{N}_h$ which intersect with a ball centred at $g\tilde{x}_0$), cover $E^f_N$ (see \eqref{22}), that is 
\[E^f_N\subset \limsup \left\lbrace  B_{\tilde{x}_0}(\eta,c_7e^{-(\rho_{g}+f(\rho_{g}))}) \;\middle|\;\begin{tabular}{@{}l@{}}
$g\in\Gamma,\eta_i\in \mathcal{J}_g$\end{tabular}\right\rbrace\] where $\eta_i$ are the centres (up to a constant at most $e^{f(\rho_g)v_{\Gamma_N}}$ many) obtained from the covering theorem. Now it follows using the fact that $f$ is Lipschitz and Lemma \ref{Borel_Cant} that if the integral in the statement of the claim converges, then $E^f_N$ has measure zero. Note that Lemma \ref{Borel_Cant} applies because the points of $\eta_i$ are uniformly radial limit points (so the volumes of small balls centred around them can be computed by the shadow lemma).

\textbf{Divergence.} 
Let $\xi^g\in S(g,R)\bigcap \partial \tilde{N}_g$. There exists a positive number $c_8>0$ large enough depending only on data (and $\epsilon$) (by $1$-thinness of $[\xi_g,z_g,gz_0]$) such that for $t>t_g:=\rho_g+c_8$ we have 
\beqn\label{25}
\tilde{\rho}(\gamma_{\xi^g}(t),\tilde{N}_g)< \epsilon/2.
\eeqn
Now if $\xi\in B(\xi^g,c_9e^{-(\rho_g+f(\rho_g))}),$ (for $c_9$ small enough depending on data and $\epsilon$) then by $1$-thinness of $[\xi,\tilde{x}_0,\xi^g]$ we have 
\beqn\label{26}
\tilde{\rho}(\gamma_\xi(t_g+f(t_g)),\gamma_{\xi^g}(t_g+f(t_g)))<\epsilon/2.
\eeqn
Therefore we have from \eqref{25} and \eqref{26} that 
\[ 
\limsup \left\lbrace B_{\tilde{x}_0}(\xi^g,c_{10}e^{-(\rho_g+f(\rho_g))})\;\middle|\;\begin{tabular}{@{}l@{}}$g\in\Gamma, \xi^g\in S(g,R)\bigcap\partial\tilde{N}_g$\end{tabular}\right\rbrace\subset E^f_N \] where $c_{10}$ depends only on data and $\epsilon$. By a $5r$-covering argument with $gz_0$ as the base point for the visual metric, we can find a disjoint collection of balls $\{B(\xi^g_i,c_{10}e^{-(\rho_g+f(\rho_g)})\}$ for varying $g$, which form a limsup set contained in $E^f_N$ to which Lemma \ref{Borel_Cant} is applied. In the divergence case, the argument in the proof of Theorem 5.1 (page 821) of \cite{HP2004}, can be used to see that the measure of $E_{f,n}$ is one when the integral diverges.
\eproof

\brem\label{rem_spiral}
In Theorem \ref{version of 4.6} we do not use the smooth structure of manifolds and the method works in the generality of CAT($-1$) spaces (see Theorem 5.3, \cite{HP3}). 
\erem

We now move on to a shrinking target problem for geodesics around totally geodesic submanifolds. This kind of theorem was first proved by Maucourant \cite{Mau} where he proved a shrinking target theorem for geodesics approximating a point in a finite volume, not necessarily compact, hyperbolic manifold. 
In particular, Theorem \ref{spiral_01} below generalises it and a theorem of Aravinda, Hersonsky and Paulin (cf. Theorem A.3 in the appendix to \cite{HP3}). Our methods also apply to cuspidal excursions, see section \ref{sec:cusp case}.\\

We need a lemma first. Recall the definition of trails of sets from \S \ref{CAT_prelim2}.
\blem[Counting cosets]\label{distribution}
Let $(X,\tilde{\rho})$ be a CAT($-1$) geodesic metric space and $\Gamma$ be a discrete group of isometries acting properly on $X$. Assume property (2) for the action of\, $\Gamma$ on $X$. Let $\Gamma_N$ be a (non-trivial) subgroup of \,$\Gamma$ acting on $X$ also with property (2) and $\tilde{N}_0$ the convex-hull of its limit set $\Lambda_{\Gamma_N}$ in $X$. Assume that $\tilde{N}_0$ is non-empty. Suppose $v_{\Gamma_N}<v_\Gamma$. Let $\tilde{x}_0\in X$. Let $\{\tilde{N}_i\}_{i\in \N\cup\{0\}}$ be the $\Gamma$-orbit of $\tilde{N}_0$. Let $z_i$ be the nearest point projections from $\tilde{x}_0$ to $\tilde{N}_i$. Let $g_i\in\Gamma$ be such that $\tilde{N}_i=g_i\tilde{N}_0$ and $$\tilde{\rho}(g_i\tilde{x}_0,z_i)=\min\{\tilde{\rho}(g'\tilde{x}_0,z_i)\mid g'\in g_i\Gamma_N\}.$$ Write $$\mathcal{N}=\{(z_i,g_i)\}_{i\in \N}\subset X\times\Gamma.$$ Then for each $h_0\in\Gamma$, $L>0$, there exists $k_0=k_0(h_0,L,X)>1, R=R(\Gamma_N,\Gamma,X)>0$ such that for all $k\in\N$, $k\geq k_{0}$,
\vspace{0.3cm}
\benum
\item We have\begin{align*}\begin{split}\#&\{(z_i,g_i) \mid z_i\in T(\tilde{x}_0,B(h_0\tilde{x}_0,L))\bigcap\left(B(\tilde{x}_0,(k+1)R)\setminus B(\tilde{x}_0,kR)\right)\}\\&\hspace{3.5cm}\approx \mu_{\tilde{x}_0}(S(\tilde{x}_0,B(h_0\tilde{x}_0,L)))e^{v_\Gamma\cdot k\cdot R}.\end{split}\end{align*} In particular, $$\#\{(z_i,g_i)\mid z_i\in B(\tilde{x}_0,(k+1)R)\setminus B(\tilde{x}_0,kR)\}\approx e^{v_\Gamma\cdot k\cdot R}.$$
\vspace{0cm}
\item Moreover, if $\Gamma$ acts cocompactly, and $\Gamma_N$ convex-cocompactly, there exists $C=C(X)>0$ such that for all $z\in B(\tilde{x}_0,(k+1)R)\setminus B(\tilde{x}_0,kR)$, there exists $(z_i,g_i)\in\mathcal{N}$ such that $$\tilde{\rho}(g_i\tilde{x}_0,z)\leq C.$$
\eenum 
\elem
The previous lemma says that the number of translates of $\tilde{N}_0$ with their `tops' or nearest points (to a fixed center, $\tilde{x}_0\in X$ here) lying in a bounded annulus centered at $\tilde{x}_0$ at distance $r$ from $\tilde{x}_0$ is $O(e^{r\cdot v_\Gamma})$. Such a statement for horospheres was used by Sullivan in the context of geodesic excursions to cusps and Diophantine approximation by parabolic fixed points, see \cite{Sul}, \cite{MP}. The second assertion of the lemma is a direct geometric consequence of the cocompactness of the action and convexity of $N$.
\bproof
Towards the proof of \textit{(2)}, write for each $k\in \N$, $$\mathcal{N}_k=\{(z_i,g_i)\in\mathcal{N}\mid z_i\in B(\tilde{x}_0,kR)\}.$$ Now fix $k\in \N$ and fix $g\in \Gamma$, such that $$g\tilde{x}_0\in B(\tilde{x}_0,(k+1)R)\setminus B(\tilde{x}_0,(k+1/2)R)$$ and there is $i\in\N$, such that $g\in g_i\Gamma_N$, where $(z_i,g_i)\in\mathcal{N}_k$. Let $z_g$ be the nearest point projection of $g\tilde{x}_0$ on $g\tilde{N}_0$. Note that in this case $\Gamma\tilde{x}_0$ is an $R_0$-net, for some $R_0>0$ (that is $\tilde{M}$ is in an $R_0$ neighbourhood of $\Gamma\tilde{x}_0$).

Let $g'\in\Gamma$ be any other such isometry, that is $g'\in g_j\Gamma_N$, $$g'\tilde{x}_0\in B(\tilde{x}_0,(k+1)R)\setminus B(\tilde{x}_0,(k+1/2)R)$$ and $(z_j,g_j)\in\mathcal{N}_k$. Define $z_{g'}$ similarly to $z_g$, that is, it is the nearest point projection of $g'\tilde{x}_0$ on $g'\tilde{N}_0$.

Write $$\xi_g=\gamma_{z_iz_g}(\infty),\;\, \xi_{g}'=\gamma_{z_jz_{g'}}(\infty).$$ Consider the geodesic triangle $[\tilde{x}_0,\xi_g,\xi_{g'}]$. Note that it follows from the CAT($-1$) inequality that for any $\epsilon>0$, a large enough $R>0$ may be chosen, so that $$\dist_{\tilde{\rho}}\,(\gamma_{\tilde{x}_0\xi_g}(kR+R/4), \gamma_{z_iz_g})\leq \epsilon,$$ and $$\dist_{\tilde{\rho}}\,(\gamma_{\tilde{x}_0\xi_{g'}}(kR+R/4), \gamma_{z_jz_{g'}})\leq \epsilon.$$ Let $w_g$ and $w_{g'}$ be the nearest point projections of $\gamma_{\tilde{x}_0\xi_g}(kR+R/4)$ on $\gamma_{z_iz_g}$ and of $\gamma_{\tilde{x}_0\xi_{g'}}(kR+R/4)$ on $\gamma_{z_jz_{g'}}$ respectively. We apply this observation, along with the fact that there exists a constant $0<C'=C'(\Gamma_N,X)$, such that $$\tilde{\rho}(w_g,w_{g'})\geq C',$$ to deduce that $$\tilde{\rho}(\gamma_{\tilde{x}_0\xi_g}(kR+R/4), \gamma_{\tilde{x}_0\xi_{g'}}(kR+R/4))\geq C'/2,$$ and thus $R$ can be chosen large enough that $$\tilde{\rho}(\gamma_{\tilde{x}_0\xi_g}(kR+R/2), \gamma_{\tilde{x}_0\xi_{g'}}(kR+R/2))\geq 100\cdot R_0.$$ From this we deduce that $$\min\{\tilde{\rho}(z_g,z_{g'}),\tilde{\rho}(g\tilde{x}_0,g'\tilde{x}_0)\}\geq 50\cdot R_0.$$ We will now use the inequality above to obtain the claim.
 
Now, let $z\in B(\tilde{x}_0,(k+1)R)\setminus B(\tilde{x}_0,(k+1/2)R)$, there exists an isometry $h\in \Gamma$, such that $$\tilde{\rho}(h\tilde{x}_0,z)\leq R_0.$$ Let $l\in\N$ be such that $h\in g_l\Gamma_N$, and $(z_l,g_l)\in\mathcal{S}$. If $(z_l,g_l)\notin\mathcal{N}_k$, the second part of the claim is verified for $z$. Consider the case $(z_l,g_l)\in\mathcal{N}_k$. Let $z_h$ be the nearest point projection of $h\tilde{x}_0$ on $h\tilde{N}_0$. There exists $C''=C''(\Gamma_N, X)$ such that $$\tilde{\rho}(h\tilde{x}_0,z_h)\leq C''.$$ Let $u_h\in X$ be a point at distance $10 R_0$ from $h\tilde{N}_0$ and $z_h$ (picked from the boundary of a tubular neighbourhood for example). Let $h'\in \Gamma$ be such that $$\tilde{\rho}(h'\tilde{x}_0,u_h)\leq R_0.$$ Then $$\tilde{\rho}(h\tilde{x}_0,h'\tilde{x}_0)\leq 12R_0,$$ and thus $$\tilde{\rho}(z,h'\tilde{x}_0)\leq 15R_0,$$ and if $m\in \N$ is such that $h'\in g_m\Gamma_N$, where $(z_m,g_m)\in\mathcal{N}$, then $(z_m,g_m)\notin\mathcal{N}_k$. This proves \textit{(2)}.

For the proof of \textit{(1)}, first write $$\mathcal{F}_m=\{g\in\Gamma\mid g\tilde{x}_0\in T(\tilde{x}_0,B(h_0\tilde{x}_0,L))\,\bigcap\,\left(B(\tilde{x}_0,(m+1)R)\setminus B(\tilde{x}_0,mR)\right)\}.$$ Then note that given $R>0$ (to be determined below) there exists $k_0=k_0(h_0, L, R)$, such that $$c_1\cdot \mu_{\tilde{x}_0}(S(\tilde{x}_0,B(h_0\tilde{x}_0,L)))e^{v_\Gamma\cdot m\cdot R}\leq \#(\mathcal{F}_{m+1}\setminus\mathcal{F}_m)\leq c_2\cdot \mu_{\tilde{x}_0}(S(\tilde{x}_0,B(h_0\tilde{x}_0,L))) e^{v_\Gamma\cdot m\cdot R},$$ for $c_1,c_2>0$ and $m\geq k_0$.

Let $k\in\N$, $k>k_0+\alpha+2$, for $\alpha\in\N$ to be fixed shortly. Let $m\in\N$ such that $k_0+\alpha+1<m<k$. For $g\in\mathcal{F}_m$, write $\mathcal{I}_g$ for the set $$\{h\in\Gamma\mid h\tilde{x}_0\in g\Gamma_N\tilde{x}_0,\, h\tilde{x}_0\in T(\tilde{x}_0,B(h_0\tilde{x}_0,L))\,\bigcap\,\left(B(\tilde{x}_0,(k+1)R)\setminus B(\tilde{x}_0,kR)\right)\}.$$ Then for $R>0$ chosen large enough, 
$$\#\mathcal{I}_g\leq c_3e^{(k-m)R\cdot v_{\Gamma_N}}.$$ Note that,
\beqnn\begin{split}\#\mathcal{F}_n\leq & \sum_{m=k_0+1}^{k-\alpha}\#(\mathcal{F}_{m+1}\setminus\mathcal{F}_m)\cdot c_3e^{(k-m)R\cdot v_{\Gamma_N}}\\&+\#(\mathcal{F}_k\setminus \{h\in\Gamma\mid h\in\mathcal{I}_g, \,\text{for some}\; g\in\mathcal{F}_{k-\alpha}\}),\end{split}\eeqnn where the second summand counts the orbits by cosets of $\Gamma_N$, whose corresponding convex hulls have tops in the annulus $B(\tilde{x}_0,kR)\setminus B(\tilde{x}_0,(k-\alpha)R)$. Then, for $R>0$ large enough, \beqnn\begin{split}&c_1\cdot \mu_{\tilde{x}_0}(S(\tilde{x}_0,B(g_0\tilde{x}_0,L)))e^{kR\cdot v_\Gamma}\\&\leq  \;c_2c_3\cdot \mu_{\tilde{x}_0}(S(\tilde{x}_0,B(g_0\tilde{x}_0,L)))\sum_{m=k_0+1}^{k-\alpha}e^{mR\cdot v_\Gamma}\cdot e^{(k-m)R\cdot v_{\Gamma_N}}\\&+ \#(\mathcal{F}_k\setminus \{h\in\Gamma\mid h\in\mathcal{I}_g, \,\text{for some}\; g\in\mathcal{F}_{k-\alpha}\}).\end{split}\eeqnn 
Then for $\alpha=\alpha(\Gamma,\Gamma_N, X)>0$ we have 

$$\#(\mathcal{F}_k\setminus \{h\in\Gamma\mid h\in\mathcal{I}_g, \,\text{for some}\; g\in\mathcal{F}_{k-\alpha}\})\geq \frac{c_1}{2}\mu_{\tilde{x}_0}(S(\tilde{x}_0,B(g_0\tilde{x}_0,L)))e^{kR\cdot v_\Gamma}.$$ This concludes the proof of \textit{(1)}.

\eproof

\bthm\label{spiral_01}
Let $M$ be a manifold of dimension $n$ with pinched negative sectional curvature $-a^2\leq k\leq -1$, for $a\geq 1$, such that the covering action of $\pi_1(M)$ on $\tilde{M}$ and $\tilde{M}$ satisfy properties (1), (2), (3) and (4). Let $N'$ be a convex submanifold in $M$ of dimension $0\leq s<n$,  such that $N'\subset N$, where $N$ is a convex submanifold of positive codimension and the action of $\pi_1(N)$ on the universal cover of $N$, $\tilde{N}\subset\tilde{M}$ satisfies property (2) and the homomorphism $\pi_1(N)\xhookrightarrow{}\pi_1(M)$ induced by inclusion is non-trivial. Let $f$ be a positive Lipschitz function. Let $x_0\in M$ be fixed. Then the set \[
  E^f_{N'} =\left\lbrace v\in SM_{x_0} \;\middle|\;
  \begin{tabular}{@{}l@{}}
    $\exists\; \text{positive times}\; t_n\to\infty\;\text{such that}$\\ $\gamma_v(t_n)\subset B(N',e^{-f(t_n)})$
   \end{tabular} 
  \right\rbrace
\] has measure zero if the integral $$\int_{1}^\infty e^{-f(t)(v_\Gamma/a-s)}dt$$ converges. It has full measure, if the integral $$\int_{1}^\infty e^{-f(t)(v_\Gamma-s)}dt$$ diverges, where $\Gamma:=\pi_1(M)$ and $\gamma_v$ is the geodesic at $x_0$ at time zero with direction $v$. If $s=n-1$, $E^f_{N'}$ has full measure. 
\ethm

\bproof  

\textbf{Convergence.}
First consider the case $0\leq s<n-1$.

We carry notation from Theorem \ref{version of 4.6}. Then there exists $0<L<\infty$ such that

$$\diam(F_{N'}\cap B(\tilde{x}_0, L))\geq\min\{1,\diam(N')\},$$ and $$\text{vol}_{\tilde{N}_0}(F_{N'}\cap B(\tilde{x}_0, L))\geq\min \{1, \text{vol}(N')\}.$$

Let $v\in E^f_{N'}$. Then there exist times $t_n\to \infty$ and isometries $g_n\in\Gamma$ such that there exist lifts $\tilde{N}_{g_n}=g_n\tilde{N}_0$  and points $x_n^v\in \tilde{N}_{g_n}$, such that $$\tilde{\gamma}_v(t_n)\in B(x_n^v,e^{-f(t_n)}).$$ 

For $g_n\in\Gamma$, consider the collection of balls $\left\lbrace B(x,\frac{e^{ -f(\rho_{g_n})}}{5})\mid x\in g_nF_N\right\rbrace$ and apply the $5r$-covering theorem to obtain a disjoint subcollection of balls $\left\lbrace B(x_j^{g_n}, \frac{e^{-f(\rho_{g_n})}}{5})\right\rbrace_j$ such that $$g_nF_N\subset \bigcup_j B(x_j^{g_n}, e^{-f(\rho_{g_n})}),$$ where the cardinality of the set of indices $j$ is (within constant multiples) of $e^{f(\rho_{g_n})s}$ (for $\rho_{g_n}$ large enough). Write $$\mathcal{J}_{g_n}=\{x^{g_n}_j\}_j$$ for the collection of centers of balls in the cover obtained above.

Note that $$\tilde{\gamma}_v(t_n)\in B(x_n^v,e^{-f(t_n)})\subset B(x_j^{g_{n,}},C_1e^{-f(\rho_{g_{n,}})}),$$ where $C_1$ depends only on the Lipschitz constant of $f$.
Let $\xi^v=\tilde{\gamma}_v(\infty)$. Let $\xi_j^{g_{n}}$ be the end point in the visual boundary of the geodesic ray starting from $\tilde{x}_0$ and passing through $x_j^{g_{n}}$. Then since the triangle $[\xi_v,\tilde{x}_0,\xi_j^{g_{n}}]$ is $a$-fat, we have $$d_{\partial \tilde{M}}(\xi^v,\xi_j^{g_{n}})\leq c_2e^{-\left(\rho_{g_{n}}+\frac{f(\rho_{g_{n}})}{a}\right)},$$ where $c_2$ is a contant depending only on the curvature bounds. Consider the functions $$Dir_R(g)=B(g\tilde{x}_0, \diam(S(g,R)))$$ and $$F(g)=\bigcup_{x^g_j\in\mathcal{J}_g} B\left(\xi^g_j, c_2e^{-\left(\rho_{g_{n}}+\frac{f(\rho_{g_{n}})}{a}\right)}\right)$$ for $g\in\Gamma$. 

Then the first part of the claim follows by Lemma \ref{Borel_Cant}: here one uses property (3) of the $\pi_1(M)$ action, by which the Patterson-Sullivan measure of balls in the limit set of $\pi(M)$, in the visual boundary $\partial \tilde{M}$ is known by Sullivan's shadow lemma, when the limit set of the group is the visual boundary, as in the hypothesis (and ergodicity ensures that the Patterson-Sullivan densities corresponding to $\pi_1(M)$ and a fixed cocompact group of isometries are both constant multiples of a Hausdorff measure).

\textbf{Divergence.} For the next part, we use terminology from Lemma \ref{distribution}. The components of the preimages of $N$ in $\tilde{M}$ are denoted $\{\tilde{N}_i\}_{i\in\N\cup\{0\}}$, where $\tilde{N}_0$ is the component passing through $\tilde{x}_0$. Recall the collection $\mathcal{N}$ from Lemma \ref{distribution}. Write $$\Gamma_{\mathcal{N}}:=\{g\in\Gamma\mid (z,g)\in\mathcal{N}, \, \text{for some}\, z\in\tilde{M}\}.$$
Given $g\in\Gamma_{\mathcal{N}}$, for a number $\kappa>1$ to be determined below, consider the collection of balls $\left\lbrace B(x,\kappa e^{ f(\rho_{g})})\mid x\in gF_{N'}\cap B(g\tilde{x}_0,L)\right\rbrace$ and apply the $5r$-covering theorem to obtain a disjoint subcollection of balls $\left\lbrace B(x_j^{g}, 5\kappa e^{-f(\rho_{g})})\right\rbrace_j$ such that $$gF_N\cap B(g\tilde{x}_0,L)\subset \bigcup_j B(x_j^{g}, 5\kappa e^{-f(\rho_{g})}),$$ where the cardinality of the set of indices $j$ is again (within constant multiples depending also on $\kappa$) of $e^{f(\rho_{g})s}$ (for $\rho_{g}$ large enough). Write $$\mathcal{J}_{g}=\{x^{g}_j\}_j$$ for the collection of centers of balls in the cover obtained above. We claim that a $\kappa\geq 1$ can be chosen such that the the shadows of the balls $B(x_i^g,e^{-f(\rho_g)})$ are mutually disjoint for $x^g_i\in\mathcal{J}_g$.

Towards the claim, we first observe that the manifold $g\tilde{N}_0$ is totally geodesic and $$\dist(B(x_l^g,e^{-f(\rho_g)}),B(x_j^g,e^{-f(\rho_g)}))\geq 5\kappa\cdot e^{-f(\rho_g)},$$ for $x^g_l\neq x^g_j$. The fact that $x_l$ and $x_j$ are within a distance bounded by a constant multiple of $L$ of the point closest to $\tilde{x}_0$ (since $g\in\Gamma_{\mathcal{N}}$) can be used to show that $$\max\{\dist_{\tilde{\rho}}(x_l,\gamma_{\tilde{x}_0x_j}), \dist{\tilde{\rho}}(x_j,\gamma_{\tilde{x}_0x_l})\}\geq c(\kappa)\cdot e^{-f(\rho_g)},$$ where $c(\kappa)$ is a positive, increasing function of $\kappa$. For this consider first the geodesic line joining $x_l$ and $x_j$ and denote it $\gamma_{x_jx_l}$. Let $z_{jl}$ be the point on $\gamma_{x_jx_l}$, closest to $\tilde{x}_0$. There are essentially two cases to consider; $x_j<z_{jl}<x_l$, (where a point to the left of the inequality comes before the point to the right along $\gamma_{x_jx_l}$), and $z_{jl}<x_j<x_l$. 

Let us consider the case $z_{jl}<x_j<x_l$. Consider the triangle $[\tilde{x}_0,x_l,y_l]$, where $y_l$ is a point at distance $\tilde{\rho}(z_{jl},x_l)$ from $z_{jl}$ along $\gamma{x_jx_l}$ on the direction opposite to that of $x_l$ and a comparison triangle $[\tilde{x}_0',x_l',y_l']$ in $\h^2_{-a^2}$.
Note that $$\tilde{\rho}(z_i,z_{jl})\leq L+\delta_M,$$ where $\delta_M$ is the constant for thinness of traingles in $\tilde{M}$ and consequently, $$\tilde{\rho}(z_{jl},x_l)\leq 2L+\delta_M\leq CL,$$ for a constant $C=(L,\delta_M)>0$.
 Let $z_{jl}'$ be the point corresponding to $z_{jl}$ in $[\tilde{x}_0',x_l',y_l']$. Since (by CBB($-a^2$) inequality) $$\rho_{-a^2}(\tilde{x}_0',z_{jl}')<\min\{\rho_{-a^2}(\tilde{x}_0',x_l'),\rho_{-a^2}(\tilde{x}_0',y_l')\},$$ the point closest to $\tilde{x}_0'$ on the geodesic line in $\h_{-a^2}$ joining $y_l'$ to $x_l'$ lies in the segment between them, let us call it $w_{jl}'$, 
Note that $$\rho_{-a^2}(x_l',w_{jl'})\leq \rho_{-a^2}(x_l',y_l')\leq CL,$$ where $C=C(L,\delta_M)>0$ and that the geodesic segment from $\tilde{x}_0'$ to $w_{jl}'$ meets the geodesic joining $x_l'$ to $y_l'$ perpendicularly. So the interior angle at $x_l'$ in the triangle $[\tilde{x}_0',x_l',y_l']$ is bounded below, depending only on $L$ (this angle monotonically tends to zero as $x_l'$ tends to infinity). Consider next the point $x_j'$ in the comparison triangle, corresponding to $x_j$, and the nearest point $u_j'$ from $x_j'$ to the side in the comparison triangle joining $\tilde{x}_0'$ to $x_l'$. By the angle lower bound at the vertex $x_l'$, there exists a constant $c=c(a,L)$, such that $$\rho_{-a^2}(x_j',u_j')\geq c(a,L)\cdot \rho_{-a^2}(x_j',x_l')=c(a,L)\cdot\tilde{\rho}(x_j,x_l).$$ 

Let $v_j$ be the point on the geodesic joining $\tilde{x}_0$ to $x_l$, nearest to $x_j$, and $v_j'$, the corresponding point in $[\tilde{x}_0',x_l',y_l']$. Then $$\rho_{-a^2}(x_j',v_j')\geq \rho_{-a^2}(x_j',u_j'),$$ and hence by the CBB($-a^2$) inequality again, $$\dist_{\tilde{\rho}}(x_j,\gamma_{\tilde{x}_0x_l})=\tilde{\rho}(x_j,v_j)\geq c(a,L)\cdot \tilde{\rho}(x_j,x_l)\geq 2 c(a,L)\cdot \kappa \cdot e^{-f(\rho_g)}.$$ 
The other case is proved arguing similarly, using the curvature lower bound, we omit the details. The claim follows by choosing $\kappa$ large enough so that $2 c(a,L)\cdot \kappa\geq 1$.

Then we have that \[\limsup\left\lbrace B(\xi_j^g,c'e^{-(\rho_g+f(\rho_g))})\;\middle|\;\begin{tabular}{@{}l@{}}$g\in\Gamma_{\mathcal{N}},j\in\mathcal{J}_g$\end{tabular}\right\rbrace\subset E^f_{N'}.\] The second part of the claim now follows again from Lemma \ref{Borel_Cant} and an argument as before which shows that the measure of $E^f_{N'}$ is in fact one, when it is positive.

Finally we consider the case that $s=n-1$. In this case the first integral does not converge, because $$v_\Gamma= v_{\tilde{M}}\leq a(n-1),$$ where $v_{\tilde{M}}$ is the critical exponent of a discrete group of isometries acting cocompactly on $M$. Indeed the first equality follows from properties (1) and (4). The second is due to results in \cite{Sinai} and \cite{Marg1} (see also \cite{Marg2}). The second integral may or may not diverge. Nevertheless, the measure is full as we explain next. By Lemmas \ref{distribution} and \ref{Borel_Cant} for a fixed $L>0$ the limsup set of shadows of balls of radius $L$ centred at orbits of $\tilde{x}_0$ by $\Gamma_{\mathcal{N}}$ has positive measure. This set is contained in $E^{f}_{N'}$ for any Lipschitz $f$ by hyperbolicity of $\tilde{M}$ and the fact that $\tilde{N}_i$ separate $M$. Thus the measure of $E^{f}_{N'}$ is positive. This completes the proof.

\eproof

\brem\label{rem_shrinking}
\benum We have the following remarks regarding Theorem \ref{spiral_01}.
\item The proof showed that for a function $f$ for which the integral diverges, the geodesic ray in almost every direction hits the target infinitely many times at an `angle' close to $\pi/2$. In fact, in the divergence part of the proof, we could estimate the sizes of the shadows of only those neighbourhoods of the preimages of $N'$, which were hit by the geodesic rays from $\tilde{x}_0$ at such angles (the corresponding translates of the fundamental domains being the ones nearest to $\tilde{x}_0$). Estimating the shadows for this collection of neighbourhoods was sufficient for applying the Borel-Cantelli lemma, because there are sufficiently many such neighbourhoods: a consequence of the non-triviality of $\pi_1(N)$ and property (2) of the $\pi_1(M)$ and $\pi_1(N)$ actions, via Lemma \ref{distribution}.
\item Since the target is hit by almost every ray, at angles close to $\pi/2$, the geodesic ray spends a time proportional to $e^{-f(t)}$ in the $e^{-f(t)}$ neighbourhood of the target, at a hitting time $t>0$. This can be checked using the lower curvature bound, with arguments similar to ones used in the proof of Theorem \ref{spiral_01}. The shrinking target phenomenon is complementary to the spiral trap phenomenon (where the time spent is much larger than the thickness of the neighbourhood), in this regard.
\item We used a volume estimate for the target. Such an estimate only requires a local Ahlfors regular measure on $N'$: a measure $m_{N'}$, and $\epsilon_{N'}>0$, such that for all $0<\epsilon<\epsilon_{N'}$, and balls $B_\epsilon$ in $N$ of radius $\epsilon$, $m_{N'}(B_\epsilon)\approx \epsilon^{\dim(N')}$. So, the manifold structure is not essential for such phenomena. 
\eenum
\erem

\subsection{Dimension estimates}\label{stron_spiral}
We begin the section by discussing the flat torus.
\bexm
Let $f(t):=\tau t$, for some $\tau>0$.
Consider the flat $n$-torus. Let $\lambda$ be the closed geodesic which is the projection of the lines $\{(k_1,\ldots,k_{n-1},t)\mid t\in\R,k_i\in\Z\}$  and consider its $\epsilon$ neighbourhood for some $\epsilon>0$. Then for any other geodesic in the torus, which is not parallel to $\lambda$ and lying in the $\epsilon$ neighbourhood, the spiral trap problem has no solutions. Indeed, for solutions to exist, $f$ has to be a bounded function. 
\eexm

\bexm
Consider the following shrinking target problem in the $3$-torus with the same function $f$ as above. Let $x_0$ be the image under the covering projection of the origin. A connected component of the preimage of a geodesic passing through $x_0$ will be of the form $t\mapsto(t,\alpha t,\beta t)$ after normalization (ignoring a set of Hausdorff dimension 1). If the above geodesic is a solution to the exponential shrinking target problem (that is with $f(t)=\tau t$ for some $\tau>0$) then it can be seen that there exist $(p_n,q_n)\in (\Z\setminus\{0\})^2$ such that $$\left|\alpha-\frac{p_n}{q_n}\right|\leq ce^{-\tau|q_n|},$$ for some absolute constant $c>0$ and all $n\in\N$. By the Jarnik-Besicovitch theorem, the possible values of $\alpha$ are a set of Hausdorff dimension zero. Thus the set of directions along which the geodesics are a solution to the exponential shrinking target problem is a set of Hausdorff dimension one (cf. Corollary \ref{full_dim}).  
\eexm

\bthm\label{loglaw1}Let $\epsilon>0$ and $\tau\geq 0$ be fixed.
Let $M$ be a manifold of pinched negative sectional curvature $-a^2\leq k\leq -1$, $a\geq 1$, such that the covering action of $\pi_1(M)$ on the convex hull of $\Lambda_{\pi_1(M)}$ in $\tilde{M}$ has property (2). Let $N$ be a compact, convex submanifold and assume $\Gamma_N:=\pi_1(N)\xhookrightarrow{}{}\pi_1(M)$ is injective. Let $x_0\in M$. Then the set 
 \[
  E^{\tau}_{N} =\left\lbrace v\in SM_{x_0} \;\middle|\;
  \begin{tabular}{@{}l@{}}
    $\exists\; \text{positive times}\; t_n \to\infty\;\text{such that}$\\ $\gamma_v(t_n,t_n+\tau t_n)\subset B(N,\epsilon)$
   \end{tabular}
  \right\rbrace
\] has $$ \text{dim}_{\mathcal{H}} (E^{\tau}_{N})= \frac{v_{\Gamma}+\tau\cdot v_{\Gamma_N}}{1+\tau},$$ where $\Gamma:=\pi_1(M)$ and $\gamma_v$ is the geodesic at $x_0$ at time zero with direction $v$.
 \ethm

\bproof
The proof follows from Theorem \ref{density} and (the proof of) Theorem \ref{version of 4.6}. Indeed, consider the function $f(t)=\tau t$. Then for both the limsup sets constructed in the proof of Theorem \ref{version of 4.6}, to approximate $E^\tau_N$ from above and below, note that Theorem \ref{density} applies with $$X=\Gamma,\, \alpha_g\approx e^{-(\rho_g+\tau\rho_g)},\, \beta_g\approx e^{v_{\Gamma_N}\cdot\tau\rho_g},\, \overline{\omega}=1.$$Here $\overline{\omega}=1$ is because of property (2). 
\eproof

We now prove a Hausdorff dimension result for the shrinking target problem. As mentioned earlier, the theorem below was known in the case that $M$ is a closed manifold and the `target' is a point due to work of Hersonsky and Paulin \cite{HP1}, see also the work of Velani \cite{Velcamb} for an analogous statement in the constant curvature case. We generalise the aforementioned results to accommodate more general manifolds $M$ and more general targets $N$.  

\bthm \label{loglaw2gen} Let $\tau\geq 0$ be fixed. Let $M$ be a manifold with pinched negative sectional curvature, $-a^2\leq k\leq -1$, $a\geq 1$, and dimension $n$, such that the covering action of $\pi_1(M)$ on $\tilde{M}$ and $\tilde{M}$ satisfies properties (2), (3) and (4). Let $N$ be a convex submanifold in $M$ of dimension $0\leq s<n$, such that the action of $\pi_1(N)$ on $\tilde{N}\subset\tilde{M}$ satisfies property (2) and the homomorphism $\pi_1(N)\xhookrightarrow{}\pi_1(M)$ induced by inclusion is non-trivial.  Then given $x_0\in M$, we have that the set 
\[
  E^{\tau}_{N} =\left\lbrace v\in SM_{x_0} \;\middle|\;
  \begin{tabular}{@{}l@{}}
    $\exists\; \text{positive times}\; t_n\to\infty\;\text{such that}$\\ $\gamma_v(t_n)\subset B(N,Ce^{-\tau t_n})\;\text{for some} \; C>0$
   \end{tabular}
  \right\rbrace
\] has Hausdorff dimension $$\frac{v_\Gamma+\tau\cdot s}{1+\tau}\leq \dim_{\mathcal{H}}(E_{N}^{\tau})\leq \frac{v_\Gamma+\tau\cdot s}{1+\tau/a},$$ where $\gamma_v$ is the geodesic at $x_0$ at time zero with direction $v$.
\ethm
\bproof
The claim again follows from Theorem \ref{density} and the proof of Theorem \ref{spiral_01}. Indeed, consider the function $f(t)=\tau t$. Then for the limsup sets constructed in Theorem \ref{spiral_01}, to approximate $E^{\tau}_N$ from above and below, the Hausdorff dimensions can be estimated by Theorem \ref{density} with $$X=\Gamma_{\mathcal{N}},\, \alpha_g\approx e^{-(\rho_g+\tau\rho_g)},\, \beta_g\approx e^{s\tau\rho_g},\, \overline{\omega}=1$$ in the case of the upper bound and with $$X=\Gamma,\, \alpha_g\approx e^{-(\rho_g+\tau\rho_g/a)},\, \beta_g\approx e^{s\tau\rho_g},\, \overline{\omega}=1,$$ for the lower bound. Here, in both cases $\overline{\omega}=1,$ by Lemma \ref{distribution}.
\eproof

\subsection{Large intersections and simultaneous spiraling}
Theorem \ref{LIP} applies to the classes of sets considered in \S \ref{stron_spiral}. We illustrate this with the following theorem. 

\bthm\label{spiral_lip}
Let $M$ be a manifold of strict-negative curvature with properties (2), (3) and (4). Let $N_i$ and $N'_i$ be countable collections of compact, totally geodesic submanifolds of $M$. Let $x_0\in M$. Let $\tau_i$ be a sequence of positive numbers and $\epsilon>0$ be given. Then the set of directions \[
  E =\left\lbrace v\in SM_{x_0} \;\middle|\;
  \begin{tabular}{@{}l@{}}
    $\exists\; \text{positive times}\; t_n^{(i)},t_n'^{(i)} \to\infty\;\text{for each $i$ such that}$\\ $\gamma_v(t_n^{(i)},t_n^{(i)}+\tau_i t_n^{(i)})\subset B(N_i,\epsilon)\;\text{and}\;\gamma_v(t_n'^{(i)})\in B(N'_i,e^{-\tau_i t_n'})$\\$\;\text{for each i}$
   \end{tabular}
  \right\rbrace
\] has $$\dim_{\mathcal{H}}(E)\geq \inf_i\min\left\lbrace\frac{v_M+\tau_i\cdot v_{N_i}}{1+\tau_i},\frac{v_M+\tau_i\dim(N_i')}{1+\tau_i}\right\rbrace.$$ If $\{\tau_i\}_i$ is bounded, the dimension lower bound is positive for non-elementary $\pi_1(M)$.
\ethm

We have the following corollary to the previous results.

\bcor\label{full_dim}
Let $M$ be a manifold of pinched negative curvature and dimension $n$; $\Gamma:=\pi_1(M)$ action on $\tilde{M}$ satisfying properties (2) and (3). Let $N$ be a compact, convex submanifold in $M$ of dimension $s$, such that $0\leq s\leq n-1$ and $\pi_1(N)\xhookrightarrow{} \pi_1(M)$ is non trivial. Let $x_0\in M$. Fix $\epsilon>0$.
\benum 
\item  If $1\leq s\leq n-1$,
 \[
  E_{N} =\left\lbrace v\in SM_{x_0} \;\middle|\;
  \begin{tabular}{@{}l@{}}
    $\exists\; \text{positive times}\; t_n \to\infty\;\text{such that}$\\ $\gamma_v(t_n,t_n+\tau t_n)\subset B(N,\epsilon)\;\text{for some}\;\tau>0$
   \end{tabular}
  \right\rbrace
\] has $$ \dim_{\mathcal{H}} (E_{N})= v_\Gamma.$$ 
\item If property (4) is also satisfied,
\[
  E'_{N} =\left\lbrace v\in SM_{x_0} \;\middle|\;
  \begin{tabular}{@{}l@{}}
    $\exists\; \text{positive times}\; t_n\to\infty\;\text{such that}$\\ $\gamma_v(t_n)\subset B(N,Ce^{-\tau t_n})\;\text{for some} \; C>0\;\text{and}\;\tau>0$
   \end{tabular}
  \right\rbrace
\] has  $$ \dim_{\mathcal{H}}(E'_{N})= v_\Gamma,$$ 
\eenum
where $\gamma_v$ is the geodesic at $x_0$ at time zero with direction $v$. Both $E_N$ and $E_N'$ have Patterson-Sullivan densities zero if $0\leq s< n-1$. If $s=n-1$, $E_N$ has again measure zero, but $E_{N}'$ has full measure, if property (1) holds (if not, the measure is still positive). 
\ecor
\bproof
The dimension results follow from Theorems \ref{loglaw1}, \ref{loglaw2gen} and \ref{spiral_lip}. The measure results follow from Theorems \ref{version of 4.6} and  \ref{spiral_01}.
\eproof

\brem\label{general} 

For any of the dimension results stated above, a smooth manifold structure is not essential; see remarks \ref{rem_spiral}
 and \ref{rem_shrinking}.\erem

\subsection{The case of cusp excursions}\label{sec:cusp case}
There has been extensive work on the shrinking target problem for cusp excursions following the work of Sullivan. Kleinbock and Margulis \cite{KMinv} generalised Sullivan's results to locally symmetric spaces of finite volume, and also gave a dynamical proof of Khintchine's theorem. See also \cites{AGP, AP} for results in the non-Archimedean setting. We would like to point out that the problem of cusp excursions of geodesics in (variable) negatively curved manifolds also falls within the framework we have developed in this paper. Namely, one can obtain both measure and dimension results for cusp excursions. In order to do so, one applies the analogue of Lemma \ref{distribution} for horospheres, known in the constant negative sectional curvature case due to the work of Sullivan and Melian and Pestana. To verify such an analogue in general, one again uses mixing for the geodesic flow with respect to the Bowen-Margulis measure. Since this result can be obtained by minor and standard modifications of the arguments in \cite{Sul} and \cite{MP}, we omit the details. This along with a minor modification of Theorem 2.1 in \cite{MP} or Theorem 3.1 in our paper yields a Jarnik-Besicovitch theorem for cusp excursions. The Borel-Cantelli statement follows from a standard modification of Sullivan's method. Theorem \ref{density} also yields a Jarn\'{i}k-Besicovitch theorem for Diophantine approximation on the Heisenberg group. Let  $X = \mathbb{H}^{n}_{\mathbb{C}}$ be complex hyperbolic space endowed with the standard Riemannian metric. Then $G = \PU(n)$ is the group of holomorphic isometries of $X$. Then $\partial X \setminus\{\infty\}$ can be identified with the Heisenberg group equipped with the Carnot-Caratheodory distance and the corresponding Hausdorff measure. This is an Ahlfors regular space. We refer the reader to \cite{HP4} for background to this problem and relevant definitions. We have that $\infty$ is a parabolic point of $\Gamma := \PU(n)\Z[\sqrt{-1}]$ and its orbit is a subset of rational points on the Heisenberg group, so   $$\{\xi\in\partial X\mid \di_{CC}(\xi,g(\infty))\leq e^{-\tau\cdot D(g(\infty))},\,\text{for infinitely many}\, g\in\Gamma\}$$
corresponds to a Diophantine approximation problem in the Heisenberg group, where $\di_{CC}$ is the Carnot-Caratheodory distance and $D(g(\infty))$ is the `depth' of the `rational' geodesic line with end points $\infty$ and $g(\infty)$. In \cite{HP4}, the authors prove a Khintchine type theorem in the Heisenberg group. Applying Theorem \ref{density}, we get a Jarn\'{i}k-Besicovitch theorem for the Carnot-Caratheodory distance. See  \cite{Zheng}, for another approach to a  Jarn\'{i}k-Besicovitch theorem on the Heisenberg group using homogeneous dynamics.

As regards the large intersection property in the case of cusp excursions, the constant negative sectional curvature case follows from the work of Falconer \cite{Fal}. However, large intersection theorems for cusp excursions in the variable negative curvature do not directly follow from the work of Falconer. However, Theorem \ref{LIP} does apply and we get the following theorem.

\bthm\label{largecat}
Let $(X,\rho)$ be a $CAT(-1)$ metric space and let $\Gamma_i$ be a collection of discrete subgroups of the isometry group of $X$ such that properties (1), (2) ,(3) and (4) hold for the action of $\Gamma_i$ on $X$. Assume further that each $\Gamma_i$ has nontrivial parabolic subgroups $P_i$. Let $\xi_i$ be corresponding fixed points in $(\partial X,\di)$. Given $\tau>0$, the set 

$$E_{\tau}=\bigcap_i\{\xi\in\partial X\mid \di(\xi,g\xi_i)\leq e^{-\tau\rho(x_0,gx_0)},\,\text{for infinitely many}\, g\in\Gamma_i\}$$ has Hausdorff dimension $$\dim_{\mathcal{H}}(E_{\tau})=\frac{v_X}{1+\tau},$$ where $v_X$ is the critical exponent of the groups $\Gamma_i$.

\ethm

In particular, the Theorem above applies to the situation where $G$ is a rank $1$ semisimple Lie group, $K$ is a maximal compact subgroup of $G$, the $\Gamma_i$ are non-uniform lattices in $G$ and the $\xi_i$ correspond to points at infinity of the finite volume quotients $\Gamma_i\backslash G/K$.

Similarly, we have a Corollary in the context of Bianchi groups. This result follows from large intersection property of Falconer in $\R^2$ and the Hausdorff content density estimate from \cite{MP} for the action (by M\"obius transformations) of the Bianchi subgroups $\Gamma_d$ on the sphere ${\Sphere}^2$ and gives a refinement of the J\'{a}rnik-Besicovitch theorem in this context. Namely, we have

\bthm\label{Bianchi}
Let $\{d_i\}_i \subset \N$ be a (possibly infinite) collection of positive square-free integers. Let $\mathcal{W}_{2, d_i}(\tau)$ denote the set of points in $\R^2$ which are  $\tau$-well approximated by the collection $\{p/q:p,q\in\Z[\sqrt{-d_i}],\; \text{Ideal}(p,q)=\Z[\sqrt{-d_i}]\}$, simultaneously for each $i$. Then $\mathcal{W}_{2, d_i}(\tau)$ is in $\mathcal{G}^{2/\tau}$. In particular, $\dim_{\mathcal{H}}(\mathcal{W}_{2, d_i}(\tau)) = 2/\tau$ and the $\mathcal{H}^{2/\tau}(\mathcal{W}_{2, d_i}(\tau)) = \infty$. Moreover,
$$ \dim_{\mathcal{H}}\left(\bigcap_{i}\mathcal{W}_{2, d_i}(\tau)\right) = 2/\tau. $$
\ethm

In the case of the Carnot-Caratheodory metric in the Heisenberg group, there is again a large intersection property. The following is a consequence of Theorem \ref{LIP}. 
\bthm
Let $\{d_i\}_i \subset \N$ be a (possibly infinite) collection of positive square-free integers. Let $\mathcal{W}_{ d_i}(\tau)$ denote the set of points in $\partial \h^n_{\mathbb{C}}\setminus \{\infty\}$ which are  $\tau$-well approximated by parabolic fixed points of $PU(n)\Z[\sqrt{-d_i}]$, simultaneously for each $i$. Then,
$$ \dim_{\mathcal{H}}\left(\bigcap_{i}\mathcal{W}_{ d_i}(\tau)\right) = \frac{2n+2}{\tau}. $$
\ethm

Finally, we note that Theorem \ref{largecat} applies to the situation when $X$ is a Bruhat-Tits tree attached to a rank $1$ semisimple algebraic group over a local field of characteristic $p$. The Diophantine problem then becomes one of approximating Laurent series by ratios of polynomials, and one can formulate the corresponding large intersection problem by considering congruence quotients, or indeed even by taking quadratic extensions (cf. \cite{GR}) as in the the Theorem above. However, we refrain from spelling out the details because the boundary of the Bruhat-Tits tree admits net measures and so Falconer's techniques apply with minor modifications. See Remark (c) in \cite{Fal}.

\end{document}